\documentclass[a4paper,11pt]{article}
\usepackage{amssymb,latexsym,amsmath,amscd,amsthm}
\usepackage{amsmath,amssymb} %per escriure mates
\usepackage{amsthm} %per configurar teoremes,definicions
\usepackage{hyperref}
\usepackage{graphicx}
\usepackage{epsfig}
\usepackage{enumerate}
\usepackage{color}
\usepackage[boxed, vlined]{algorithm2e} %option noend, lined , noline
\usepackage{lscape}
\usepackage{changes}
\usepackage{appendix}

\definecolor{pink}{rgb}{1,0,1}
 
\definechangesauthor[color=pink]{js}
\definechangesauthor[color=blue]{mrn}
\definechangesauthor[color=red]{mt}

\theoremstyle{definition}
\newtheorem{defi}{Definition}[section]

\theoremstyle{plain}
\newtheorem{lema}[defi]{Lemma}
\newtheorem*{lema*}{Lemma}

\newtheorem{thm}[defi]{Theorem}

\newtheorem{assumption}[defi]{Assumption}

\newtheorem{prop}[defi]{Proposition}
\newtheorem{cor}[defi]{Corollary}
\newtheorem{rk}[defi]{Remark}

\newtheorem{teo-def}[defi]{Theorem/Definition}
\newtheorem{alg}[defi]{Algorithm}

\newcommand{\ZZ}{\mathbb{Z}}
\newcommand{\CC}{\mathbb{C}}

\newcommand{\RR}{\mathbb{R}}

\newcommand{\dd}{\mathrm{diag}}

\newcommand{\Arg}[1]{\mathrm{Arg}({#1})}
\newcommand{\Log}{\mathrm{Log}}
\newcommand{\Ll}[2]{\mathrm{Log}_{#1}(M)}%{\mathrm{Log}(#1;#2)}

\newcommand{\Disc}{\Delta(x)}

\newenvironment{proofThmintro}{\noindent {\textit{Proof of Theorem 1.1.}}}{\hfill $\square$ \vspace{3mm}}

\DeclareMathOperator{\trace}{tr}
\DeclareMathOperator{\Real}{Re}
\DeclareMathOperator{\Imag}{Im}
\newif\ifprivate
 \privatetrue
\def\???{\ifprivate {\bf {???}} \marginpar{{\Huge {\bf ?}}}
\else \fi}

\definecolor{garnet}{RGB}{210,15,30}
\newcommand{\jr}[1]{\textcolor{garnet}{#1}}

\title{The embedding problem for Markov matrices}
\author{Marta Casanellas, Jes\'us Fern\'andez-S\'anchez, Jordi Roca-Lacostena}
\date{}

\begin{document}

\maketitle

\begin{abstract}
%Characterizing when a general Markov process has a a continuous-time realization is a hard problem, which remains open for $k\geq 4$ states. In practice, this problem reduces to deciding when a given $k\times k$ Markov matrix can be written as the exponential of some rate matrix (Markov generator). This question is known in literature as the embedding problem \cite{Elfving}.
%

%In this paper, we address this and related questions and obtain results in two different lines. First, for general $k$, we give a bound for the number of Markov generators in terms of the spectrum of the Markov matrix and based on this, we establish a criterion for deciding whether a generic Markov matrix is embeddable. We also propose an algorithm that lists all Markov generators.
%
%Then, motivated and inspired by recent results on substitution models of DNA, we focus in the case $k=4$ and completely solve the embedding problem for any Markov matrix. The solution in this case is more elegant as the embeddability is given in terms of a single condition.
%
%
%

Characterizing whether a Markov process of discrete random variables has a homogeneous continuous-time realization is a hard problem. In practice, this problem reduces to deciding when a given Markov matrix can be written as the exponential of some rate matrix (a Markov generator). This is an old question known in the literature as the \textit{embedding problem} \cite{Elfving}, which has been only solved for matrices of size $2\times 2$ or $3\times 3$.
In this paper, we address this problem and related questions and obtain results in two different lines. First, for matrices of any size, we give a bound on the number of Markov generators in terms of the spectrum of the Markov matrix. Based on this, we establish a criterion for deciding whether a generic (distinct eigenvalues) Markov matrix  is embeddable and propose an algorithm that lists all its Markov generators.
Then, motivated and inspired by recent results on substitution models of DNA, we focus in the $4\times 4$ case and completely solve the embedding problem for any Markov matrix. The solution in this case is more concise as the embeddability is given in terms of a single condition.

\end{abstract}


\section{Introduction}
Markov matrices are used to describe changes between the states of two discrete random variables in a Markov process. As the entries of Markov matrices (or transition matrices) represent the conditional probabilities of substitution between states, Markov matrices have non-negative entries and rows summing to one. Among them, embeddable matrices are those that are consistent with a homogeneous continuous-time Markov process, so that changes occur at a constant rate over time and time is thought as a continuous concept. 
{The instantaneous rates of substitution are usually displayed as the entries of real matrices with non-negative off-diagonal entries and rows summing to zero, so-called \emph{rate matrices}. In the homogeneous continuous-time setting, the transition matrices of a Markov process can be computed as $M(t)=e^{Qt}$, where $Q$ is the rate matrix ruling the process and $t\geq 0$ accounts for the time elapsed in the process. In this case, $M$ is said to be \emph{embeddable}. Equivalently, a Markov matrix $M$ is embeddable if it can be written as the exponential of a rate matrix $Q$, $M = e^Q$ (with no reference to time $t$). Any rate matrix $Q$ satisfying $M = e^Q$ is called a \emph{Markov generator} of $M$.}

Almost one century ago, Elfving \cite{Elfving}  formulated the problem of deciding which Markov matrices are embeddable, \emph{the embedding problem}. Solving the embedding problem results in giving necessary and sufficient conditions for a Markov matrix $M$  to be the exponential of a rate matrix $Q$, $M=e^Q$. Although the question is quite theoretical, it has practical consequences and, as such, it may appear in every applied field where discrete and continuous-time Markov processes are considered. For instance, in economic sciences \cite{Israel,Geweke}, in social sciences \cite{Singer} and in evolutionary biology \cite{verbyla,ChenJia}, the embedding problem is crucial for deciding whether a Markov process can be modeled as a homogeneous continuous-time process or not.

Although the embedding problem is solved for $2\times 2$ and $3 \times 3$ matrices \cite{Kingman,Cuthbert73,Johansen,Carette}, it has remained open for larger matrices so far. Some partial results on the necessary conditions for a Markov matrix to be embeddable were given in the second part of the twentieth century \cite{Runnenburg,Kingman,Cuthbert72}. Moreover, there exist sufficient and necessary conditions on the embeddability of Markov matrices with \emph{different and real} eigenvalues. This is a consequence of a result due to Culver \cite{Culver} and characterizes embeddability of this type of matrices in terms of the principal logarithm, see Corollary \ref{cor:DistinctRealEmbed}.
%states that a Markov matrix with positive (NO CAL?) different eigenvalues is embeddable if and only if its principal logarithm is a rate matrix (Markov matrices with negative eigenvalues of multiplicity one are never embedabble).
%
There are also some inequalities that need to be satisfied by the determinant or the entries of the matrix in order to be embeddable \cite{Goodman,Fuglede}.
%MENCIONAR MES TREBALLS PER MATRIUS GRANS? (antics Cuthver...),
At the same time, there is a discrete version of the embedding problem, which consists on deciding when a Markov matrix can be written as a certain power of another Markov matrix (see \cite{Singer,Guerry2013,Guerry2019} for instance).
%This discrete version has received some attention in a number of papers (see \cite{Singer,Guerry2013,Guerry2019} for instance).

A related issue is deciding whether there is a unique Markov generator for a given embeddable Markov matrix. Note that each Markov generator provides a different embedding of the Markov matrix into a homogeneous continous-time Markov process. We refer to this question as the \emph{rate identifiability problem}. It is well known that for diagonally dominant embeddable matrices, the number of Markov generators reduces to one  \cite[Thm 4]{Cuthbert72}. The same happens if the matrix is close to the identity; for example, if either $|| M-I ||<0.5$ or $\det(M)>0.5$  \cite{Israel}. However, the situation becomes really complicated as the determinant of the matrix decreases. The first example of a Markov matrix with more than one Markov generator was given in \cite{Speakman}, and further examples were provided in \cite{Cuthbert73,Israel}. In all these examples, however, the principal logarithm happens to be a rate matrix.% and so, it is one of the Markov generators of the Markov matrix.

%DIR EL QUE FEM PER MATRIUS GRANS. TECNIQUES QUE USEM:
In this paper we provide a solution to the embedding problem for Markov matrices of any size with pairwise \emph{different} eigenvalues (not necessarily real), see Theorem \ref{thm:nxnCriterion}.
This situation covers a dense open subset of the space of Markov matrices, so it solves the embedding problem almost completely  (the set of matrices with repeated eigenvalues has measure zero within the whole space of matrices).
For such matrices, we bound the number of Markov generators in terms of the real and imaginary parts of the eigenvalues and establish a criterion for deciding whether a Markov matrix with different eigenvalues is embeddable. Based on this criterion, we provide an algorithm that gives \emph{all} Markov generators for Markov matrices with different eigenvalues (Algorithm \ref{Alg:nxn}). We also give an improvement in the  bounds on the determinant mentioned above, see Corollary \ref{cor:DetBound}. The main techniques are the description of the complex logarithms of a matrix (see \cite{Gantmacher}) and a careful study of the complex region where the eigenvalues of a rate matrix lie (Section 3).

In addition to these results, we completely solve the embedding problem for $4\times 4$ Markov matrices (with \emph{repeated or different} eigenvalues). The solution to the embedding problem provided this case (see Section 5) is much more satisfactory because we are able to characterize embeddability by checking a single condition (and not looking at a list of possible Markov generators).
We have devoted special attention to $4\times 4$ matrices not only because it was the first case that remained still open, but also because our original approach and motivation arises from the field of phylogenetics, where Markov matrices rule the substitution of nucleotides in the  evolution of DNA molecules.  In the last years, new results and advances concerning the embedding problem have appeared in this field, providing deep insight and illustrative examples of the complexity of the general situation, see \cite{ChenJia, JJ, Baake, K2}.  The present work builds on some previous contributions by the authors in this setting.

For $4\times 4$ Markov matrices $M$ with different eigenvalues (real or not) we prove that the embeddability can be checked directly by looking at the principal logarithm $\Log(M)$ together with  a basis of eigenvectors:
%matrices our main result  is Corollary \ref{cor:ConjugateEmbeddingCriterion}, which together with Lemma \ref{lem_cas1} solves the embedding problem for Markov matrices $M$ with different eigenvalues. The result states that the embeddability of $M$ can be checked directly by looking at its principal logarithm and its eigenvectors:
 %comparing three values $\mathcal{L}$, $\mathcal{U}$ and $\mathcal{N}$ which can be computed directly from $M$.
\begin{thm}\label{thm_intro} Let $M=P\dd(1,\lambda_1,\lambda_2,\lambda_3)P^{-1}$ be a $4\times 4$ Markov matrix with $\lambda_1\in \RR_{>0}$, $\lambda_2\in \CC$, $\lambda_3 \in \CC$ pairwise different. If $\lambda_2, \lambda_3 \notin \RR$, define $V= P \;\dd(0,0,2\pi i ,-2\pi i)  \; P^{-1}$,
$$\mathcal{L}:= \displaystyle \max_{(i,j):\ i\neq j,\  V_{i,j} > 0} \left\lceil -\frac{\Log(M)_{i,j}}{V_{i,j}}  \right\rceil,  \qquad \mathcal{U}:= \displaystyle \min_{(i,j):\ i\neq j, \  V_{i,j} < 0} \left\lfloor -\frac{\Log(M)_{i,j}}{V_{i,j}}
  \right\rfloor ,$$
and define $V=0$, $\mathcal{L}=\mathcal{U}=0$ if all eigenvalues are real. Set   \[\mathcal{N}:= \{(i,j): i\neq j, \  V_{i,j}=0 \text{ and  } \Log(M)_{i,j}<0\}.\]
Then, $M$ is embeddable if and only if $\mathcal{N} = \emptyset$, $\mathcal{L} \leq \mathcal{U}$ and $\lambda_i \notin \RR_{\leq 0}$ for $i=1,2,3$. In this case, the Markov generators of $M$ are {the matrices $\Log(M)+ k\, V$ where $k\in\mathbb{Z}$ satisfies $\mathcal{L}\leq k \leq \mathcal{U}$.} 
\end{thm}

As a byproduct we give an algorithm that outputs \emph{all} possible Markov generators for such a matrix. Apart form this general case of matrices with different eigenvalues, we also study all other cases and we give an embeddability criterion for each (see Section 5.1, cases I, II, III, IV, and Section 5.2). The case of diagonalizable matrices with \emph{two} real repeated eigenvalues (Case III) turns out to be much more involved; still we are able to provide necessary and sufficient conditions for the embeddability in terms of eigenvalues and eigenvectors, and to propose an algorithm that checks whether a Markov matrix in this case is embeddable (Cor. \ref{cor:RepeatedEmbedabilitIFF}, Alg. \ref{Alg2}).

The outline of the paper is as follows. In Section \ref{sec:Prelim} we state with precision the embedding problem and recall some known results needed in the sequel.
% In Section \ref{sec:Counterexample} we show that there is a set of positive measure containing embeddable matrices whose principal logarithm is not a rate matrix.
Section \ref{sec:bounds} is devoted to bounding the real and the imaginary part of the eigenvalues of any rate matrix (Lemma \ref{lema:EigenBound}).
%As an easy consequence, we improve the bound on the determinant given in  \cite{Israel} that guarantees $\mathrm{Log}(M)$ as the unique possible Markov generator \blue{M: Aixo ho podem treure perque ja ho hem dit abans}.
%a technical result by \cite{Israel} to check embeddability in terms of the principal logarithm.
These bounds are used in Section \ref{sec:NoRepeated} in order to provide a sufficient and necessary condition for an $n\times n$ Markov matrix with pairwise different eigenvalues to be embeddable.
%In practice, this result (Theorem \ref{thm:nxnCriterion}) solves the embedding problem for generic Markov matrices.
In the same section, we also give the algorithm that outputs all Markov generators of such matrices.
%we state a sufficient and necessary condition for embeddability for those $4\times4$ matrices that do not diagonalize or have at least three equal eigenvalues. Moreover we provide an algorithm that given a $4\times4$ Markov matrices with two repeated eigenvalues (Algorithm \ref{Alg2}) decides whether it is embeddable or not and if it has a finite or infinite number of Markov generators solving also the rate-identifiability problem.
We devote Section \ref{sec:4x4} to $4\times 4$ matrices, studying their embeddability with full detail by splitting them into all possible Jordan canonical forms. The proof of Theorem \ref{thm_intro} is also given there.
%Finally in section \ref{sec:Conclusions} we discuss the impact of the results of the paper.
%
In the last section of the paper, Section \ref{sec:RateId}, we summarize the results on  the {rate identifiability} for embeddable $4\times 4$ matrices (see Table \ref{table_aux}).
%
%
%Appendix \ref{sec:Counterexample} contains the proofs of those results of the previous sections that involve different techniques; namely, some topological arguments for Theorem \ref{Thm:Appendix} and Proposition \ref{Prop:AllkM}.
%
%
% and handling polytopes for the implementation of Algorithm \ref{Alg2}.
Appendix  \ref{sec:Case3} is devoted to details concerning the implementation of Algorithm \ref{Alg2}.

\section{Preliminaries}\label{sec:Prelim}

%!TEX root = main.tex"`

In this section we recall some definitions and relevant facts about the embedding problem of Markov matrices.

\begin{defi}
A real square matrix $M$ is a \textit{Markov matrix} if its entries are non-negative and all its rows sum to 1. A real square matrix $Q$ is a \textit{rate matrix} if its off-diagonal entries are non-negative and its rows sum to 0. A Markov matrix $M$ is \textit{embeddable} if there is a rate matrix $Q$ such that $M=e^Q$; in this case we say that $Q$ is a \textit{Markov generator} for $M$. Embeddable Markov matrices are also sometimes referred to as matrices that have a \emph{continuous realization} \cite{steelbook}. The \textit{embedding problem} \cite{Elfving} consists on deciding whether a given Markov matrix is embeddable or not, in other words determine  which Markov matrices can be embedded into the multiplicative semigroup $\Big( \{e^{Qt}:t\geq 0 \} , \cdot \Big)$ for some rate matrix $Q$.
\end{defi}

%\begin{defi}
%An embeddable Markov matrix $M$ has identifiable rates if there exists a unique rate matrix $Q$ such that $M=e^Q$.  The \textit{rate identifiability problem} consists on deciding whether a given Markov matrix has identifiable rates or not.
%\end{defi}

The following notation will be used throughout the paper.
%\begin{notat}\rm
$Id$ denotes the identity matrix of order $n$. We write $GL_n(\mathbb{K})$ for the space of $n\times n$ invertible matrices with entries in $\mathbb{K}= \RR$ or $\CC$.  For $\lambda \in \mathbb{C}\setminus \{0\}$, we use the notation $\log_k(\lambda)$ to denote the $k$-th branch of the logarithm of $\lambda$, that is, $\log_k(\lambda)=\log|\lambda|+ (\text{Arg}(\lambda) + 2\pi k) i$ where $\text{Arg}(\lambda) \in (-\pi, \pi]$ is the \emph{principal argument} of $\lambda$. For ease of reading the principal logarithm $\log_0(\lambda)$ will be denoted as $\log(\lambda)$. Given a square matrix $M$, we denote by $\sigma(M)$ the set of all its eigenvalues and by $Comm^*(M)$ the \emph{commutant} of $M$, that is, the set of invertible complex matrices that commute with $M$.
%\end{notat}

\begin{rk}\label{rk_comm}
\rm If $D$ is a diagonal matrix,  $D=\dd(\overbrace{\lambda_1,\dots,\lambda_1}^{m_1},\overbrace{\lambda_2,\dots,\lambda_2}^{m_2},\dots,\overbrace{\lambda_l,\dots, \lambda_l}^{m_l})$ with $\lambda_i\neq \lambda_j$, then $Comm^*(D)$ consists of all the block-diagonal matrices whose blocks are taken from the corresponding $GL_{m_i}(\CC)$.
In particular, the commutant of $D$ does not depend on the particular values of the entries $\lambda_i$.
If $m_1=m_2=\dots=m_l=1$  then $Comm^*(D)$ is the set of invertible diagonal matrices.
\end{rk}

 If $M$ is diagonalizable, the following result describes all possible \emph{logarithms} of $M$  (that is, all the solutions $Q$ to the equation $M=e^Q$).

\begin{thm}[  {\cite[\S VIII.8]{Gantmacher}}] \label{thm:CharOfLog}
Given a non-singular matrix $M$ with an eigendecomposition $ P \; \dd(\lambda_1,\lambda_2,\dots,\lambda_n) \; P^{-1}$, where $\lambda_i \in \CC$, $i=1,\dots, n$, and  $P\in GL_n(\CC)$, the following are equivalent:
\begin{enumerate}[i)]
\item $Q$ is a solution to the equation $M=e^Q$,

\item \label{eq:Qki}
 $Q =P \; A\; \dd\big(\log_{k_1}(\lambda_1),\log_{k_2}(\lambda_2),\dots,\log_{k_n}(\lambda_n) \big) \; A^{-1} \; P^{-1}$ for some $k_1,k_2,\dots,k_n\in \ZZ$ and some $A\in Comm^* (\dd(\lambda_1,\lambda_2,\dots,\lambda_n))$.

\end{enumerate}
\end{thm}

%From now on, we will denote the preceding matrix by $\log_{k_1,k_2,\dots,k_n}^{A}(M)$.

\begin{rk} \label{rk:IfDifNoCommutator}
\rm With respect to the previous result, we want to point out the following.
\begin{enumerate}
 \item[(i)] If $u$ is an eigenvector of $Q$ with eigenvalue $a$, then $u$ is also an eigenvector of $M=e^Q$ with eigenvalue $e^a$. The converse is not true in general.
 %The eigenvalues of any logarithm $Q$ of $M$ are logarithms of the eigenvalues of $M$; the eigenvectors of $Q$ are also eigenvectors of $M$ (the converse is not true in general).
 \item[(ii)]  If $Comm^*(\dd(\lambda_1,\dots,\lambda_n))=Comm^*(\dd\big(\log_{k_1}(\lambda_1),\dots,\log_{k_n}(\lambda_n) \big))$ then the description of the logarithms is slightly simpler, as every logarithm can be written as
\[Q= P \; \dd\big(\log_{k_1}(\lambda_1),\dots,\log_{k_n}(\lambda_n) \big) \; P^{-1}.\]Moreover, in this case, $M$ and $Q$ have the same eigenvectors.
%Some relevant cases where the equality of the commutants above does hold are: i) if
This occurs, for example, when all the eigenvalues of $M$ are pairwise distinct or also when $k_1=k_2=\dots=k_n$.% In these cases, $M$ and $Q$ have the same eigenvectors.
\end{enumerate}

%This will be reflected in the notation by writing $log_{k,k,\dots,k} (M)$ instead of $\log_{k,k,\dots,k}^P(M)$.
%Furthermore, if $k_i=k_j$ for all $\lambda_i=\lambda_j$ then $A$ also commutes with $\dd\big(\log_{k1}(\lambda_1),\log_{k2}(\lambda_2),\dots,\log_{kn}(\lambda_n) \big)$ thus $\log_{k1,k2,\dots,kn}^{P}(M)$ does not depend on $P$ other than for fixing the ordering of the eigenvalues.
%Note that in all these cases, $M$ and these logarithms of $M$ have the same eigenvectors.

%For ease of reading we will use the notation $\log_{k1,k2,\dots,kn}(M)$ assuming that the eigenvalues are ordered as follows: first the real eigenvalues ordered by decreasing norm, and later the pairs of conjugated complex eigenvalues ordered by decreasing norm and each pair ordered by its imaginary part.
\end{rk}

%\private{Comentar que estem extenent la def de Log(M) per agafar vaps negatius}
%\begin{rk}
%\tp{Note that this definition of the principal logarithm extends the usual definition (see \cite{Higham}) that requires that the matrix $M$ has no negative eigenvalues. Since all the arguments in this paper work independently of this restriction, we have decided to generalize the idea of principal logarithm for the sake of simplicity and clarity in the statements. }
%\end{rk}

{The well-known formula $\det e^Q=e^{\mathrm{tr}\,  Q}$ implies that the determinant of every embeddable matrix is a positive real number. Hence, throughout this paper  we implicitely assume that \emph{all Markov matrices are non-singular} and \emph{have positive determinant}. Note that this is not a restriction of the original problem, but a necessary condition for a Markov matrix to be embeddable.}

Since both Markov and rate matrices have only real entries, the study about the existence of real logarithms of real matrices by \cite{Culver} is relevant for solving the embedding problem. The following proposition is a direct consequence of that work.

\begin{prop}[see {\cite[Thm. 1]{Culver}}]\label{prop:Culver}

Let $M$ be a real square matrix. Then, there exists a real
logarithm of $M$ if and only if $\det(M)>0$ and
each Jordan block of $M$ associated with a negative eigenvalue  occurs an even number of times.
\end{prop}

\vspace{2mm}
{In \cite[Theorem 2]{Culver}, Culver also proved that matrices with pairwise distinct positive eigenvalues have only one real logarithm, namely, the \emph{principal logarithm}}:

\begin{defi}\label{def:principa}
% Given two matrices $M,Q$ such that $M=e^Q$ we say that $Q$ \textit{is a logarithm of }$M$.
%We say that $Q$ \textit{is a primary logarithm of} $M$ if $Q=\log_{k,k,\dots,k}(M)$ for any $P$ that diagonalizes $M$ .
The \textit{principal logarithm of} $M$, which will be denoted by $\Log(M)$, is the only logarithm whose eigenvalues are the principal logarithm of the eigenvalues of $M$ (see \cite[Thm 1.31]{Higham}). In particular, if $M$ is diagonalizable then
\[\Log(M)=P\;\dd(\log(\lambda_1),\dots,\log(\lambda_n))\;P^{-1}.\]%
If $M$ is a Markov matrix, then its principal logarithm $\Log(M)$ has row sums equal to 0 (although it may not be a real matrix).
\end{defi}

\begin{rk}\rm
Note that the above definition of the principal logarithm (Definition \ref{def:principa}) extends the usual definition (e.g. see \cite[pp. 20]{Higham}), which requires that the matrix $M$ has no negative eigenvalues. This is required in order to use the spectral resolution of the logarithm function. In this paper, however, we  mainly deal with diagonalizable matrices, for which the principal logarithm can be defined directly by taking the principal argument of the eigenvalues.  The only non-diagonalizable Markov matrices that we deal with are $4\times 4$ (see Section 5.2) which, according to Prop. \ref{prop:Culver}, have no negative eigenvalues if they have a real logarithm.
\end{rk}
%\private{EL QUE DIC ES CERT, PERO NO ES IMMEDIAT. SI TINGUES VAP NEGATIU FORA DELS BLOCS DE JORDAN, AMB EL VAP 1, ELS  2 VAPS NEGATIUS I EL BLOC DE JORDAN, JA ENS PASSEM DE DIMENSIO. CAL EXPLICAR-HO? Marta: NO, NO CAL EXPLICAR-HO, FEM REFERENCIA A AQUELLA SECCIO I ALLA JA ES VEU. HE CANVIAT REDACTAT. TOT I AIXI EN LA DEF 2.5 HEM D'ASSEGURAR-NOS QUE ENCARA QUE HI HAGI VAPS NEGATIUS EN UNA MATRIU DIAGONALITZABLE, NOMES HI HA UN UNIC LOG AMB ARGUMENT EN (-pi,pi]; ES CERT?}

  %(this was already shown by \cite{Singer}).
{As a byproduct of the results explained above,} we get the following embeddability criterion for Markov matrices with pairwise distinct real eigenvalues  in terms of its principal logarithm.

\begin{cor}\label{cor:DistinctRealEmbed}
Let $M$ be a Markov matrix with pairwise distinct real eigenvalues. Then:\begin{enumerate}[i)]
\item If $M$ has a non-positive eigenvalue, then $M$ is not embeddable.
\item If $M$ has no negative eigenvalues, $M$ is embeddable if and only if $\Log(M)$ is a rate matrix.
\end{enumerate}
\end{cor}

\section{Bounds on the eigenvalues of rate matrices}\label{sec:bounds}
It is well known that the eigenvalues of a Markov matrix have modulus smaller than or equal to one \cite[\S 8.4]{Meyer}. Here we bound the real and the imaginary part of the complex eigenvalues of rate matrices. To this end, if $Q$ is an $n \times n$ rate matrix with $n\geq 3$ and $\lambda \in \sigma(Q)$ is a non-real eigenvalue, we define
$$b_n(\lambda):=\min \left \{\sqrt{2 \trace(Q) \Real(\lambda) - (\Real(\lambda))^2}, -\frac{\Real(\lambda)}{\tan(\pi/n)} \right \},$$
$$B_n:=\min \left \{ -\frac{\sqrt{3}}{2} \trace(Q), -\frac{\trace(Q)}{2\tan(\pi/n)}\right \}.$$

The following technical result is used in the next section and is also useful to improve a result of \cite{Israel} (see Corollary \ref{cor:DetBound}).

\begin{lema}\label{lema:EigenBound}
Let $Q$ be a $n\times n$ rate matrix. Then for any eigenvalue $\lambda \in \sigma (Q)$ we have
\begin{enumerate}
 \item[i)] $\Real(\lambda)\leq 0$. Moreover, if $\lambda \not \in \RR$ then $\frac{\trace(Q)}{2} \leq \Real(\lambda)\leq 0$.
 \item[ii)] $|\Imag(\lambda)| \leq b_n(\lambda) \leq B_n$ if $\lambda \notin \RR.$
 %\item[ii)] $|\Imag(\lambda)| \leq \sqrt{2 \trace(Q) \Real(\lambda) - (\Real(\lambda))^2} \leq -\frac{\sqrt{3}}{2} \trace(Q)$ for all $\lambda \not \in \RR$.
 %\item[iii)] $|\Imag(\lambda)| \leq-\frac{\Real(\lambda)}{\tan(\pi/n)} \leq -\frac{\trace(Q)}{2\tan(\pi/n)} $ for all $\lambda \not \in \RR$.
\end{enumerate}
{Moreover, the bound on $|\Imag(\lambda)|$ given by $b_n(\lambda)$ is tight for $n\geq 3$.}
\end{lema}

\begin{proof}
%\begin{enumerate}
%\item[i)]
$i)$ If $Q$ is a rate matrix then $e^Q$ is a Markov matrix. In particular, the eigenvalues of $Q$ are logarithms of the eigenvalues of a Markov matrix. Since the modulus of the eigenvalues of a Markov matrix is bounded by $1$, we get $\Real(\lambda)\leq 0$ for any $\lambda \in \sigma(Q)$. Moreover, as non-real eigenvalues of $Q$ appear in conjugate pairs, we have that $$\trace(Q) = \sum _{\lambda\in \sigma(Q)}\lambda = \sum _{\lambda\in \sigma(Q)\cap \RR}\lambda + \sum _{\lambda\in \sigma(Q)\setminus \RR} \Real(\lambda).$$
Therefore, if $\lambda\notin \RR$, then $\Real(\lambda)$ appears twice in this expression, and so
$\Real(\lambda)\geq\trace(Q)/2$.
%This proof is based on the proof of Thm 5.1 c) by \cite{Israel}.
%\item[ii)]

\vspace*{2mm}

$ii)$ We prove first that, for any non-real eigenvalue $\lambda \in \sigma(Q)$, we have
\begin{equation}\label{eq_bound1}
|\Imag(\lambda)| \leq \sqrt{2 \trace(Q) \Real(\lambda) - (\Real(\lambda))^2} \leq -\frac{\sqrt{3}}{2} \trace(Q).
\end{equation}

 Let us take $r = -\trace(Q)$. Since $Q$ is a rate matrix we get that $\widetilde{Q}=Q+rId_n$ is a matrix with non-negative entries whose rows sum to $r$.
    %It follows from the Gershgorin Theorem (see Theorem B.1 in \cite{Higham}) that for any $\widetilde{\lambda} \in \sigma(\widetilde{Q})$ there is a row $i$ of $\widetilde{Q}$ such that $|\widetilde{\lambda}-\widetilde{q_{ii}}|\leq \sum_{j\neq i}\ |\widetilde{q_{ij}}|$. However, using that $\widetilde{Q}$ has non-negative entries, we get that for any row $i$ we have:
 %$$ |\widetilde{\lambda}|-\widetilde{q_{ii}}=|\widetilde{\lambda}|-|\widetilde{q_{ii}}|\leq |\widetilde{\lambda}-\widetilde{q_{ii}}|\leq \sum_{j\neq i}\ |\widetilde{q_{ij}}| =  \sum_{j\neq i}\ \widetilde{q_{ij}}.$$
Then any eigenvalue $\widetilde{\lambda}\in \sigma(\widetilde{Q})$ has modulus smaller than or equal to $r$  (see \cite[\S 8.3]{Meyer}).
%$|\widetilde{\lambda}| \leq \sum \widetilde{q_{ij}} = r$.
Now, if $\lambda$ is an eigenvalue of $Q$ we have that $\lambda+r \in \sigma(\widetilde{Q})$. Therefore, $(\Real(\lambda)+r)^2+\Imag(\lambda)^2= | \lambda+r|^2\leq r^2$ and we obtain \begin{equation}\label{eq:ImBound}
|\Imag(\lambda)| \leq \sqrt{r^2- (\Real(\lambda)+r)^2}= \sqrt{2 \; \Real(\lambda) \; \trace(Q)  - \Real(\lambda) ^2}.
\end{equation}
The second inequality in \eqref{eq_bound1} follows by using $0 \geq \Real(\lambda) \geq  \trace(Q)/2$ in (\ref{eq:ImBound}). %

%\item[iii)]
We prove now that
\begin{equation}\label{eq_bound2}
  |\Imag(\lambda)| \leq-\frac{\Real(\lambda)}{\tan(\pi/n)} \leq -\frac{\trace(Q)}{2\tan(\pi/n)}
\end{equation}
for any non-real eigenvalue $\lambda$ of $Q$.
 If $n<3$ then $Q$ has no complex eigenvalue, because $0$ is an eigenvalue of any rate matrix and complex eigenvalues of real matrices appear in conjugate pairs. If $n\geq3$, the first theorem in \cite{Runnenburg} claims that the principal argument of any eigenvalue $\lambda \in \sigma(Q)$ is bounded as \begin{equation}\label{eq:EigenOfRate}
\left(\frac{1}{2}+\frac{1}{n}\right)\pi\leq |\Arg\lambda|.
\end{equation}
Then the first inequality in \eqref{eq_bound2} is obtained by using that $\Imag(\lambda)=\tan(\Arg{\lambda}) \Real(\lambda)$, $\Real(\lambda)\leq 0$, and that $|\tan|$ restricted to $(-\pi,-\pi/2-\pi /n] \cup [\pi/2+\pi/n,\pi]$ attains its maximum at
\begin{equation}\label{eq:Tangent}
\left| \tan \left(-\frac{\pi}{2}-\frac{\pi}{n} \right) \right| =\left| \tan \left(\frac{\pi}{2}+\frac{\pi}{n} \right) \right|= \frac{1}{\tan(\pi /n)} > 0.
\end{equation}
%\tan \left(\left(\frac{3}{2}-\frac{1}{n}\right)\pi \right)=
The second inequality follows by using $\Real(\lambda) \geq  \trace(Q)/2.$%

From the first inequality in both \eqref{eq_bound1} and \eqref{eq_bound2}, we have  $|\Imag(\lambda)| \leq b_n(\lambda)$. On the other hand, the inequality $b_n(\lambda) \leq B_n$ follows from the definition of $b_n(\lambda)$ and the second inequalities in \eqref{eq_bound1} and \eqref{eq_bound2}. This concludes the proof of (ii).

Next, we show that the bound on $|\Imag(\lambda)|$ given by $b_n(\lambda)$ is tight for $n\geq 3$. As shown in \cite{Runnenburg} (Theorem of page 537), for each  $n\geq 3$, there is an $n\times n$ rate matrix $Q$ with at least one eigenvalue satisfying \eqref{eq:EigenOfRate} with equality. Such a matrix is given by
\begin{equation*} %\label{eq:QTight}
q_{i j} = \begin{cases}
- \alpha & \text{if } i=j,\\
 \alpha & \text{if } i\equiv j-1 \mod n,\\
0 & \text{otherwise.}\\
\end{cases}
\end{equation*}
where $\alpha$ is an arbitrary positive number. 
It can be seen that this matrix $Q$ has at least a non-real eigenvalue $\lambda$ satisfying
$|\mathrm{Im}(\lambda)|=\left| \Real(\lambda)\tan \left( \frac{\pi}{2} +\frac{\pi}{n} \right) \right|.$ Using that $\Real(\lambda)<0$ together with \eqref{eq:Tangent}, we get  $-\frac{\Real(\lambda)}{\tan(\pi/n)} = |\mathrm{Im}(\lambda)| \leq b_n(\lambda)$.
%\leq -\frac{\Real(\lambda)}{\tan(\pi/n)}$, which concludes the proof.
The definition of $b_n(\lambda)$ implies that $|\mathrm{Im}(\lambda)|=b_n(\lambda)$, so the inequality is tight. 
%
%Now, from \eqref{eq_bound1} and \eqref{eq_bound2} we obtain $|\Imag(\lambda)| \leq b_n(\lambda) \leq B_n.$
%%\end{enumerate}
%%%%%%%%%%%%%%%%5
%%\begin{equation}\label{eq:betterEigenBound}
%
%Finally, as we have seen that for any non-real eigenvalue $\lambda$ of $Q$, $\trace(Q) \geq 2 \Real(\lambda)$, we get
%$${\small \sqrt{2 \trace(Q) \Real(\lambda) - (\Real(\lambda))^2} \geq \sqrt{4 \Real(\lambda)^2 - (\Real(\lambda))^2} = -\sqrt{3} \Real(\lambda)} %= -\frac{\Real(\lambda)}{\tan(\pi/6)}. $$
%
%The final claim follows by observing that $\tan(\pi/n)\geq \tan(\pi/6)$ for any $n \in \{3,4,5,6\}$.
%
\end{proof}

\begin{figure}
	\centering
	\includegraphics[width=\linewidth]{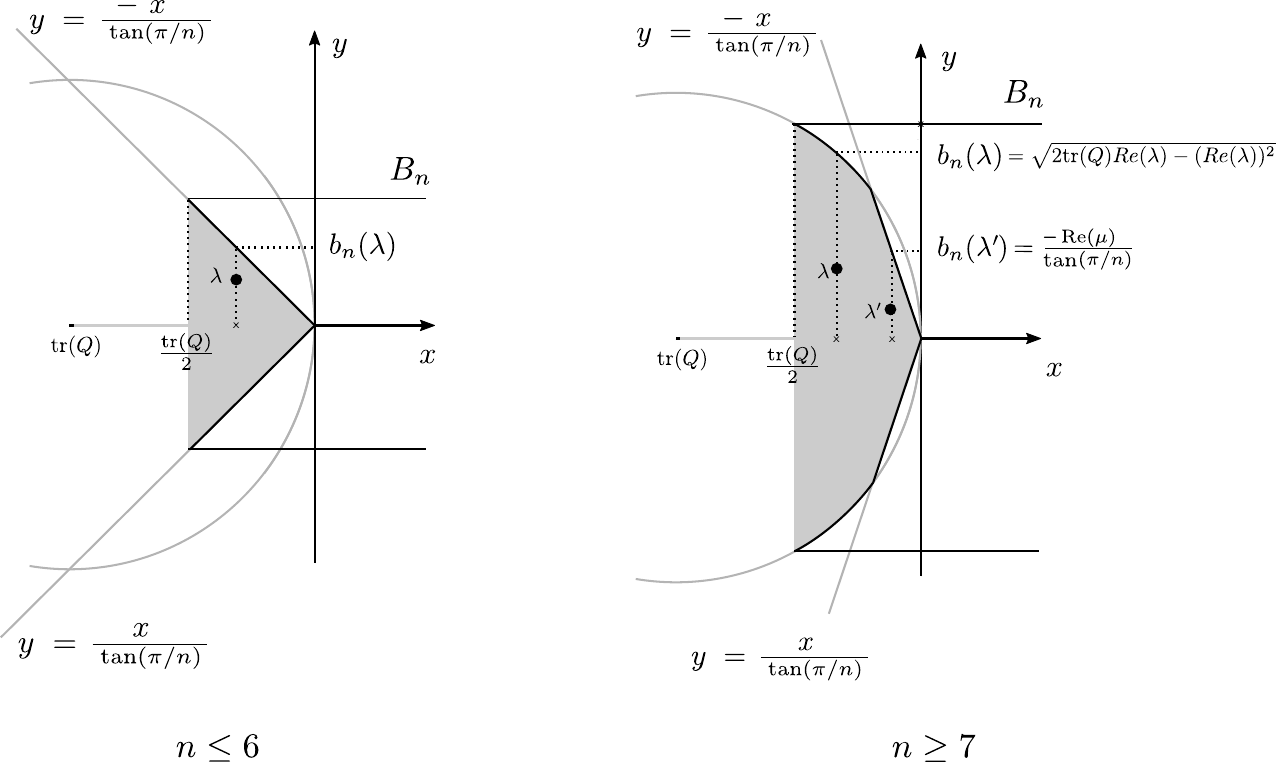}
	\caption{\label{Im:cercle}
	%The gray color shows the area where the eigenvalues of a rate matrix $Q$ lie (\ref{lema:EigenBound} ii) ). }
	We represent in dark the complex region containing the eigenvalues of $n\times n$ rate matrices (Lemma \ref{lema:EigenBound}) according to the two cases described in Remark \ref{rk:BnValue}. %
	For $n\leq 6$, we have that $B_n=-\frac{tr(Q)}{2 \tan(\pi/n)}$ and $b_n(\lambda)=-\frac{\mathrm{Re}(\lambda)}{\tan(\pi/n)}$; while for $n\geq 7$, $B_n=-\frac{\sqrt{3}}{2}\mathrm{tr}(Q)$ and the expression for $b_n(\lambda)$ depends on the value of $\lambda$ as illustrated by the eigenvalues $\lambda$ and $\lambda'$ in the figure. %
	All these bounds are represented in black in the figure.}
\end{figure}

\begin{rk}\label{rk:BnValue}
\rm
%Note that for $n>6$ the value of $b_n(\lambda)$ is not decided \tp{determined?} only by $n$ but also by the value of $\lambda$.
Note that % if $3\leq n\leq 6$, then $b_n(\lambda)=-\frac{\Real(\lambda)}{\tan(\pi/n)}$.
for $n\geq 2$, we have $\tan(\pi/n) > \tan(\pi/(n+1))$. As for $n=6$ $\tan(\pi/n)$ is $1/\sqrt{3}$, we obtain that
\begin{equation}\label{eq:bn}
B_n= \begin{cases}

-\frac{\trace(Q)}{2\tan(\pi/n)} & \text{ if } n=3,4,5,6\\

-\frac{\sqrt{3}}{2} \trace(Q) & \text{ if } n\geq 6 \, .
\end{cases}
\end{equation}
 Although for $n>6$ the value in the minimum defining $b_n(\lambda)$ depends on $\lambda$, for $n\leq 6$ we have $b_n (\lambda) = -\frac{\Real(\lambda)}{\tan(\pi/n)}$. Indeed, if $n\in\{3,4,5,6\}$, then
$${\footnotesize
	-\frac{\Real(\lambda)}{\tan(\pi/n)} \leq
	{- \Real(\lambda)\sqrt{3}  =
	\sqrt{4 \Real(\lambda)^2 - \Real(\lambda)^2} }\leq
	\sqrt{2 \trace(Q) \Real(\lambda) - (\Real(\lambda))^2}
}.$$
Figure \ref{Im:cercle} illustrates both the remark above and Lemma \ref{lema:EigenBound}.
\end{rk}

The following result improves the bound given in \cite[Theorem 5.1]{Israel}, which states that a Markov matrix $M$ with pairwise distinct eigenvalues and $\det(M)> e^{-\pi}$ is embeddable if and only if $\Log(M)$ is a rate matrix. We are able to relax the hypothesis on the determinant and avoid the condition of distinct eigenvalues.

\begin{cor}\label{cor:DetBound}
Let $M$ be a $n\times n$ Markov matrix with $\det(M) > \min \left \{e^{-\frac{2 \pi}{\sqrt{3}}},e^{-2 \pi \tan(\pi/n)} \right \}$. Then, the only possible Markov generator for $M$ is $\Log(M)$. In particular, $M$ is embeddable if and only if $\Log(M)$ is a rate matrix.
%Then, the unique possible Markov generator of $M$ is $\Log(M)$. In particular, $M$ is embeddable if and only if $\Log(M)$ is a rate matrix.
\end{cor}

\begin{proof}
Let $Q$ be a Markov generator for $M$. By hypothesis, $\trace(Q)=\log(\det(M))$ is strictly greater than $\min\{{-\frac{2 \pi}{\sqrt{3}}},{-2 \pi \tan(\pi/n)\}}.$ Therefore, using Lemma \ref{lema:EigenBound}(ii), we have $|\Imag(\lambda)| \leq   B_n <\pi$ for all $\lambda \in \sigma(Q)$. Hence, $Q$ is the principal logarithm of $M$.
\end{proof}

\begin{rk}\label{rk:DetBound}
\rm
As in Remark \ref{rk:BnValue}, we have that  $e^{-2 \pi \tan(\pi/n)} \leq e^{-\frac{2 \pi}{\sqrt{3}}}  $ for $n=3,4,5,6$ and   $e^{-2 \pi \tan(\pi/n)} \geq e^{-\frac{2 \pi}{\sqrt{3}}}  $ for $n\geq 6$.

% We believe that the condition on the determinant could be relaxed. As far as we are aware, for $n=4$, the largest determinant of an embeddable matrix with a generator different than $\Log(M)$ is $e^{-4\pi}$ (see \cite{K2},  Rmk. 4.6).
% % Therefore we cannot expect the bound in Cor. \ref{cor:DetBound} to be improved below $e^{-4\pi}$ (for $n=4$).
% Moreover, note that the bound in Corollary \ref{cor:DetBound} arises from $B_n$ in Lemma \ref{lema:EigenBound}. Hence, a more relaxed  hypothesis depending not only on the determinant of $M$ but also on its eigenvalues could be obtained by using $\max_{\lambda \in \sigma(M)} b_n( \log|\lambda|)$ instead of $B_n$.

% In the case of statement iii) this is done in  Corollary \ref{cor:EigenGlobalBound}.
 %However, these new bounds and the bounds provided in Corollaries \ref{cor:DetBound}  and \ref{cor:TangentDetBound} coincide if all the eigenvalues of $M$ excluding a conjugate pair  have modulus equal to 1. If the bound is sharp, the limit case belongs to this family of matrices and the conjugate pair of eigenvalues must be a negative eigenvalue with multiplicity 2.
\end{rk}

Table \ref{tab:DetBounds} gives numerical values for the bound in Corollary \ref{cor:DetBound}.% with other previously known bounds for several sizes of Markov matrices.

\begin{table}[h!]
  \centering
\begin{tabular}{|c | c c c c |}
\hline
Size of $M$  & $n=3$ & $n=4$ & $n=5$ & $n\geq6$ \\
 \hline
%\cite{Cuthbert73} &   0.043214 & 0.043214 & 0.043214 & 0.043214 & $0.043214$ \\
% \cite{Israel} & 0.5 & 0.5 & 0.5 & 0.5 & 0.5 \\
%$e^{-2\pi/\sqrt{3}}$ & 0.026580 & 0.026580 & 0.026580 & \textbf{0.026580} & \textbf{0.026580} \\
Bound on $\det(M)$ & 0.000019 &  0.001867 &  0.010410 &  0.026580 \\
\hline
\end{tabular}
\caption{\label{tab:DetBounds} Lower bounds on the determinant (rounded to the 6th decimal) that allow the
characterization of the embeddability in terms of the principal logarithm of an $n\times n$ Markov matrix according to Corollary \ref{cor:DetBound} and Remark \ref{rk:DetBound}. %For $n\geq 6$ then bound is $e^{-2\pi/\sqrt{3}}= 0.026580$.
Previously known bounds were $e^{-\pi}=0.043214$ and $0.5$ for all $n\in \mathbb{N}$ (see \cite{Cuthbert73} and \cite{Israel},  respectively), which have been improved by Corollary \ref{cor:DetBound}.}
\end{table}

\section{Embeddability of Markov matrices with (non-real) distinct eigenvalues}\label{sec:NoRepeated}
%!TEX root = main.tex"`

In this section we deal with Markov matrices with pairwise distinct eigenvalues. It is known that the embeddability of these matrices is determined by the principal logarithm if all the eigenvalues are \textit{real} (see Corollary \ref{cor:DistinctRealEmbed}). However,  this may not be the case if there is a non-real eigenvalue \cite{Birkhauser}. %Indeed,
% the forthcoming Proposition \ref{prop:nxnNecessary} shows that
{By virtue of Proposition \ref{prop:Culver}, we know that every (real) matrix $M$ with distinct eigenvalues has some real logarithm; by \cite[Theorem 2 and Corollary]{Culver}, we also know that if $M$ has complex eigenvalues, then there are countable infinitely many logarithms. In this section, we give a precise description of those real logarithms with rows summing to zero (Proposition \ref{prop:nxnNecessary}) and
we show that only a finite subset of them have non-negative off-diagonal entries (Theorem \ref{thm:nxnCriterion})}. In this way we are able to design an algorithm that returns the Markov generators of a Markov matrix with distinct eigenvalues (real or not), see Algorithm \ref{Alg:nxn}.
%Moreover we provide a criterion to determine the embeddability of a Markov matrix with distinct eigenvalues and some of them being non-real.  To this end, we first bound the number of logarithms.
%As matrices with non-positive determinant have no real logarithm \cite{Culver}, we will work with matrices with positive determinant.

{As pointed out in the Preliminaries section, it is well known that a necessary condition for a Markov matrix to be embeddable is to have positive determinant. If we assume that all eigenvalues are distinct, then this implies that the real eigenvalues lie in the interval $(0,1]$.}
% It is well known that if such a Markov matrix  has some non-positive eigenvalue then it has no Markov generator (Proposition \ref{prop:Culver}) and if it is singular it does not even have a logarithm. On the other hand, all eigenvalues of a Markov matrix have modulus bounded from above by one.  Thus, we are only interested in  Markov matrices whose real eigenvalues lie in the interval $(0,1]$.
{Throughout this section we consider matrices satisfying the following assumption:}

\begin{assumption}\label{Assump}
{ We assume that $M$ is a (non-singular) diagonalizable $n\times n$ Markov matrix with pairwise distinct eigenvalues and whose real eigenvalues lie in the interval $(0,1]$. Without loss of generality, we can write} %$M$ is a Markov matrix with an eigendecomposition of the form
\begin{equation*}%\label{eqM}
M=P\; \dd \big( 1,\lambda_1,\dots,\lambda_t, \mu_1,\overline{\mu_1},\dots,\mu_s,\overline{\mu_s} \big) \; P^{-1}
\end{equation*}
for some $P\in GL_n(\CC)$, $\lambda_i \in (0,1)$ with $i=1,\dots,t$, and $\mu_j \in \{z\in \CC: \Imag(z) > 0\}$ with $j=1,\dots,s$, all of them pairwise distinct.
\end{assumption}

\begin{defi}
%Given $M=P\; \dd \big( 1,\lambda_1,\dots,\lambda_t, \mu_1,\overline{\mu_1},\dots,\mu_s,\overline{\mu_s} \big) \; P^{-1}$ an $n\times n$ Markov matrix with pairwise distinct eigenvalues, $\lambda_i \in (0,1)$ for $i=1,\dots,t$ and $\mu_j \in \{z\in \CC: \Imag(z) > 0\}$ for $j=1,\dots,s$, $P\in GL_n(\CC)$,
Given a Markov matrix $M$ as in Assumption \ref{Assump}, for each $(k_1,\dots,k_s)\in \ZZ^s$ we define the following matrix:
\small
\begin{eqnarray*}
\Log_{k_1,\dots,k_s}(M):=P\; \dd \Big( 0, \log(\lambda_1),\dots,\log(\lambda_t), \log_{k_1}(\mu_1),\overline{\log_{k_1}(\mu_1)},\dots,\log_{k_s}(\mu_s),\overline{\log_{k_s}(\mu_s)} \Big) \; P^{-1}.
\end{eqnarray*}
\normalsize
Note that $\Log_{0,\dots,0}(M)$ is the principal logarithm of $M$, $\Log(M)$.
\end{defi}

The next result claims that these are all the real logarithms of the matrix $M$.
% uses de definition above to characterize  all the real logarithms with rows summing to $0$ of any given $n\times n$ Markov matrix with no repeated eigenvalues.

\begin{prop}\label{prop:nxnNecessary}
%Let $M=P\; \dd \big( 1,\lambda_1,\dots,\lambda_t, \mu_1,\overline{\mu_1},\dots,\mu_s,\overline{\mu_s} \big) \; P^{-1}$ be an $n\times n$ Markov matrix with pairwise distinct eigenvalues, $\lambda_i \in (0,1)$ for $i=1,\dots,t$ and $\mu_j \in \{z\in \CC: \Imag(z) > 0\}$ for $j=1,\dots,s$, $P\in GL_n(\CC)$.
Let $M$ be a Markov matrix as in Assumption \ref{Assump}.
Then, a matrix $Q$ with rows summing to $0$ is a real logarithm of $M$ if and only if $Q=\Log_{k_1,\dots,k_s}(M)$ for some $k_1,\dots,k_s \in \ZZ$.\\

\end{prop}

%\private{This proof should be much simpler.}

\begin{proof}

 We know that the first column of $P$ is an eigenvector of $M$ with eigenvalue $1$. Since the rows of $M$ sum to one and it has no repeated eigenvalue we can assume without loss of generality that it is the eigenvector $(1,1,\dots,1)$. In addition, we know that for $m=1,\dots,s$ the $(t+2m)$-th and $(t+2m+1)$-th columns of $P$ can be chosen to be conjugates because $M$ is real.

\begin{enumerate}

\item[$\Leftarrow$)] Since $\overline{\log_{k}(\mu)}= \log_{-k}(\overline{\mu})$ it follows from Theorem \ref{thm:CharOfLog} that $\Log_{k_1,\dots,k_s}(M)$ is a logarithm of $M$ for any $k_1,\dots,k_s \in \ZZ$. Note that the rows of $Q$ sum to $0$ because the  first column of $P$ is the eigenvector $(1,1,\dots,1)$  and its corresponding eigenvalue is $0$. Moreover, the non-real eigenvalues of $Q$ appear in conjugate pairs and the corresponding  eigenvectors appearing as column-vectors in $P$ are also {a conjugate pair}, thus $Q$ is real.

\item[$\Rightarrow$)] Let  $Q$ be a real logarithm of $M$ with rows summing to $0$.
 Since $M$ has pairwise distinct eigenvalues so does $Q$. Moreover, $Q$ {is diagonalizable} through $P$ (see Remark \ref{rk:IfDifNoCommutator} (ii)). Hence, it follows from Theorem \ref{thm:CharOfLog} that:
 \begin{eqnarray*}
 Q & = & P\; \dd \Big( \log_{k_0}(1), \log_{k_1}(\lambda_1),\dots,\log_{k_t}(\lambda_t), \dots \\
 \hspace*{10mm} & &  \dots, \log_{k_{t+1}}(\mu_1),\log_{k_{t+2}}(\overline{\mu_1}),\dots,\log_{k_{t+2s-1}}(\mu_s),\log_{k_{t+2s}}(\overline{\mu_s}) \Big) \; P^{-1}
 \end{eqnarray*}
 Since the rows of $Q$ sum to $0$ we get that $k_0 =0$. Since $Q$ is real and has no repeated eigenvalues it follows that $k_1=k_2=\dots=k_t=0$ and that its non-real eigenvalues appear in conjugate pairs. Hence, $\log_{k_{t+2m-1}}(\overline{\mu_m}) =\overline{\log_{k_{t+2m}}(\mu_m)}$.

\end{enumerate}
\end{proof}

\begin{rk}
 \rm
%  Note that the proposition above holds even if $s=0$ (that is, all the eigenvalues of $M$ are real) and in this case, $\Log_{k_1,\dots,k_s}(M)$ is the principal logarithm of $M$.
When all the eigenvalues of $M$ are real (that is, $s=0$), the proposition above claims that the only real logarithm with rows summing to 0 is the principal logarithm.
\end{rk}

From the proposition above and Lemma \ref{lema:EigenBound} we get that any Markov matrix with pairwise distinct eigenvalues has a finite number of Markov generators. Hence, its embeddability can be determined by checking whether a finite family of well-defined matrices contains a rate matrix or not, as stated in the next result. In order to simplify the notation, for a given Markov matrix $M$ and for any $z\in \CC$ we define

$$\beta_n(z):=\min \left \{\sqrt{2\log(\det(M))\log|z|-\log^2|z| }, -\frac{\log|z|}{\tan(\pi/n)} \right \}
.$$

If $Q$ is a Markov generator of $M$ and $\log_k(z)$ is an eigenvalue of $Q$ then $\beta_n(z) =  b_n(\log_k(z))$. Hence, according to Lemma \ref{lema:EigenBound} we have $\beta_n(z)=-\frac{\log|z|}{\tan(\pi/n)}$ for $n=3,4,5,6$ .

\begin{thm}\label{thm:nxnCriterion}
If $M$ is a Markov matrix as in Assumption \ref{Assump}, then
%Let $M$ and  $s$ be defined as in Proposition \ref{prop:nxnNecessary}. Then:
 \begin{enumerate}
\item[i)]$M$ is embeddable if and only if $\Log_{k_1,\dots,k_s}(M)$ is a rate matrix for some $(k_1,\dots,k_s)\in \ZZ^s$ satisfying $
\left\lceil \frac{-\Arg{\mu_j}- \beta_n(\mu_j)}{2\pi} \right\rceil\leq k_j \leq \left\lfloor \frac{-\Arg{\mu_j} +\beta_n(\mu_j)}{2\pi} \right\rfloor$ for $j=1,\dots,s$.

\item[ii)] $M$ has at most $ \left \lfloor 1 -\frac{\sqrt{3}\log(\det(M)) }{2 \pi} \right \rfloor ^s$ Markov generators if $n\geq 6$, at most $\left \lfloor 1 -\frac{\log(\det(M))}{2\pi \tan(\pi/n)} \right \rfloor ^s$ if $n=3,4,5$ and at most one if $n\leq 2$.  \label{stat:genBound}
\end{enumerate}
\end{thm}

\begin{proof}

If $Q$ is a logarithm of $M$, then $\log(\det(M))=\trace(Q)$. Hence, since $\Real(\log_{k}(\mu))= \log|\mu|$ we have that $\beta_n(\mu) = b_n(\log_{k}(\mu))$  for any $k \in \ZZ$ and any $\mu \in \sigma(M)$.
\begin{enumerate}
\item[i)]Let $Q$ be a Markov generator for $M$. According to Proposition \ref{prop:nxnNecessary} there exist $k_1,\dots,k_s\in \ZZ$ such that $Q=\Log_{k_1,\dots,k_s}(M)$. Now, by Lemma  \ref{lema:EigenBound} we have $|\Imag\big(\log_{k_j}(\mu_j)\big)|  \leq  \beta_n(\mu_j)$.
 We get the asserted bounds by using that $ |\Imag\big(\log_{k_j}(\mu_j)\big)| = | \Arg{\mu_j} + 2\pi k_j| $.
%  : $$\left\lceil \frac{-Arg(\mu_j)- \beta_j}{2\pi} \right\rceil\leq k_j \leq \left\lfloor \frac{-Arg(\mu_) +\beta_j}{2\pi} \right\rfloor.$$

\item[ii)] If $n<3$ then $M$ has only real eigenvalues and hence its only possible Markov generator is $\Log(M)$. For other values of $n$, it follows from the first statement that if $\Log_{k_1,\dots,k_s}(M)$ is a Markov generator then $k_j$ lies in an interval of length $\frac{2\beta_n(\mu_j)}{2\pi}$. Since $k_j\in \ZZ$ for all $j$ we get that $M$ has at most $\prod_j \left\lfloor  1+ \frac{2\beta_n(\mu_j)}{2\pi} \right\rfloor$ generators.

The statement follows by using Lemma \ref{lema:EigenBound} and Remark \ref{rk:BnValue} to get
$$\beta_n(\mu) =
%b_n(\log_{k}(\mu))\leq B_n=
 \begin{cases}

-\frac{\log(\det(M))}{2\tan(\pi/n)}   & \text{ if } n=3,4,5,6\\

-\frac{\sqrt{3}}{2} \log(\det(M)) & \text{ if } n\geq 6.

\end{cases}$$

  %Using that $\trace(Q)= \log(\det(M))$, the statement follows from \ref{lema:EigenBound} ii).

%On the other hand, using that $\log(\det(M))=\trace(Q)$ and Lemma \ref{lema:EigenBound}, we get $\beta(\mu_j) \leq -\frac{\sqrt{3}}{2}\log(\det(M))$, which concludes the proof.
\end{enumerate}
\end{proof}

%\begin{comment}
%\begin{cor}\label{cor:genBound2}
%Let $M$ and $s$ be defined as in Proposition \ref{prop:nxnNecessary}. Then, $M$ has at most $\left \lfloor 1 -\frac{\sqrt{3}\log(\det(M)) }{2 \pi} \right \rfloor ^s$ Markov generators.

%\end{cor}

%\begin{proof}
%It follows from Theorem \ref{thm:nxnCriterion} that if $\Log_{k_1,\dots,k_s}(M)$ is a Markov generator then $k_j$ lies in an interval of length $\frac{2\beta(\mu_j)}{2\pi}$. Since $k_j\in \ZZ$ for all $j$ we get that $M$ has at most $\prod_j \left\lfloor  1+ \frac{2\beta(\mu_j)}{2\pi} \right\rfloor$ generators. On the other hand, it follows from Lemma \ref{lema:EigenBound} that $\beta(\mu_j) \leq -\frac{\sqrt{3}}{2}\log(\det(M))$, which concludes the proof.
%\end{proof}
%\end{comment}

\begin{rk}
\rm
As shown in the proof of Theorem \ref{thm:nxnCriterion} $ii)$, the number of Markov generators of $M$ is also bounded by $\prod_j \left\lfloor  1+ \frac{2\beta_n(\mu_j)}{2\pi} \right\rfloor$. Although this bound improves  those in Theorem \ref{thm:nxnCriterion} $ii)$, we do not know if it is sharp or not and we preferred to give a bound depending on $\log(\det(M))$ because this quantity might be related to the expected number of substitutions of the Markov process ruled by  $M$ (see \cite{BarryHartigan} for further details on this in the context of phylogenetics).% Both bounds might coincide if $M$ has  at most one pair of conjugated eigenvalues ($s=1$).
\end{rk}

To close this section, we present an algorithm which determines the embeddability of a Markov matrix with pairwise distinct eigenvalues and returns all its Markov generators.\\

\begin{alg}[Markov generators for $n\times n$ matrices with distinct eigenvalues]\label{Alg:nxn}\

\begin{algorithm}[H]%\label{Alg:nxn}\
\SetKwData{Left}{left}\SetKwData{This}{this}\SetKwData{Up}{up}
\SetKwFunction{Union}{Union}\SetKwFunction{FindCompress}{FindCompress}
\SetKwInOut{Input}{input}\SetKwInOut{Output}{output}

\Input{$M$, an $n\times n$ Markov matrix with no repeated eigenvalues.}
\Output{All its Markov generators if $M$ is embeddable, an empty list otherwise.}
\BlankLine
generators=[\ ]\

compute eigenvalues of $M$\

\If{M has a negative or zero eigenvalue}{
\textbf{return} ``M not embeddable''\

\textbf{exit}}

\Else{
	s = $\frac{\# \text{non-real eigenvalues}}{2}$\
	
	\If{$s >0$ (i.e. $M$ has a non-real eigenvalue)}
	{ 	
		\For{$j=1,\dots,s$}
		{set $l_j =  \left\lceil \frac{-\Arg{\mu_j} -\beta_n(\mu_j)}{2\pi} \right\rceil$  and  $u_j =  \left\lfloor \frac{-\Arg{\mu_j} +\beta_n(\mu_j)}{2\pi} \right\rfloor$ }
  		\For{$k_1=l_1,\dots,u_1$}{
	 	 $\ddots$\
	 	
		\For{$k_s=l_s,\dots,u_s$}
				{compute $\Log_{k_1,\dots,k_s}(M)$\
				
   					\If{$\Log_{k_1,\dots,k_s}(M)$ is a rate matrix}
   					{\textbf{add} $\Log_{k_1,\dots,k_s}(M)$ to generators}
   				}
   			
   		}
 	}
	\Else
	{
		\If{$\Log(M)$ is a rate matrix}
		{\textbf{add} $\Log(M)$ to generators}
	}

	\If{generators=[\ ]}
	{\textbf{return} ``M not embeddable''}
	\Else{\textbf{return} generators}
}
%\caption{Markov generators for matrices with different eigenvalues.}
\end{algorithm}
\end{alg}

\begin{rk}\rm
As stated in Corollary \ref{cor:DetBound}, if $M$ has a Markov generator different than $\Log(M)$, then $M$ has a small determinant and some eigenvalues of $M$ are close to $0$. In this case there might be numerical issues in the implementation of the algorithm.
%because the real part of the eigenvalues of $\Log_{k_1,\dots,k_s}(M)$  is computed as $\log_{k_j}|\mu_j|$ and this is a large value which might have large errors depending on  the accuracy used to compute $\mu_j$.
\end{rk}

\section{Embeddability of $4\times 4$ Markov matrices}\label{sec:4x4}
%!TEX root = main.tex"`

In this section we study the embedding problem for all $4\times 4$ Markov matrices. In this case, we can be more precise than in Theorem \ref{thm:nxnCriterion} and, for matrices with distinct eigenvalues,  we manage to give a criterion for the embeddability in terms of the eigenvectors, see Corollary \ref{cor:ConjugateEmbeddingCriterion}. We are also able to deal with repeated eigenvalues so that the results of this section include all possible  $4\times 4 $ Markov matrices.

%As claimed in the Introduction, we already know that the embeddability of a Markov matrix is not always determined by the principal logarithm. In the $4\times 4$ case, we can prove that the set of embeddable Markov matrices whose $\Log(M)$ is not a Markov generator is not a subset of zero measure; on the contrary, it is a set of full dimension. Since the techniques needed to prove this are different to those used so far, the proof is delayed to the Appendix.

%\begin{thm}[Appendix, Theorem \ref{thm:openset}]\label{Thm:Appendix}
%There is a non-empty Euclidean open set of $4\times 4 $ \emph{embeddable} Markov matrices whose principal logarithm is not a rate matrix.
%\end{thm}

%\begin{proof}
%See Theorem \ref{thm:openset} in Appendix.
%\end{proof}

%This section is divided into two. First, we deal with the case of diagonalizable matrices with full detail. Then, we will focus on the case of remaining matrices.

\subsection{Embeddability of diagonalizable $4\times4$ Markov matrices}

%From Gershgorin Theorem we know that the modulus of all eigenvalues of a Markov matrix is bounded by one. Moreover, by Perron-Frobenius Theorem we know that if $M$ is positive (i.e. it has no null entries) then the only eigenvalue with modulus $1$ is $1$ itself and it has algebraic multiplicity $1$.
{Using Proposition \ref{prop:Culver} and the fact that the modulus of the eigenvalues of a Markov matrix is bounded by 1, we are able to} enumerate all possible diagonal forms of a diagonalizable $4\times4$ Markov matrix with real logarithms (up to ordering the eigenvalues):\\

\begin{lema}\label{lemma:DiagCases}
Let $M$ be a diagonalizable $4\times 4$ Markov matrix. If $M$ admits a real logarithm then its diagonal form lies necessarily in one of the following cases (up to ordering the eigenvalues):\\

%\begin{table}[h]
\begin{tabular}{l c l}
\rm \textbf{Case I} & $\dd(1,\lambda_1,\lambda_2,\lambda_3)$ & with $\lambda_1,\lambda_2,\lambda_3\in(0,1]$  pairwise distinct.\\
\rm \textbf{Case II} & $\dd(1,\lambda,\mu,\bar{\mu})$ & with $\lambda\in(0,1]$, $\mu,\bar{\mu} \in \CC\setminus\RR$.\\
\rm \textbf{Case III} & $\dd(1,\lambda,\mu,\mu)$ & with $\lambda\in (0,1]$, $\mu\in [-1,1)$, $ \mu\neq 0$, $\mu \neq \lambda$.\\

\rm  \textbf{Case IV} & $\dd(1,\lambda,\lambda,\lambda)$ &  with $\lambda\in (0,1]$.\\
\end{tabular}
\end{lema}

\begin{proof}
%First note that $M$ is non-singular (Proposition \ref{prop:Culver}).
%Recall that   $|\lambda|\leq 1$ for any $\lambda \in \sigma(M)$ and $(1,1,1,1)^t$ is an eigenvector of $M$ with eigenvalue $1$  (Proposition \ref{prop:PF}).
%
Since $M$ is diagonalizable and $1$ is an eigenvalue, we can write $\dd(1,\lambda_1,\lambda_2,\lambda_3)$ for the diagonal form.
If $M$ has a negative eigenvalue, it must have multiplicity $2$ by Proposition \ref{prop:Culver}. Thus, $M$ has at most one negative eigenvalue.
On the other hand, since $M$ is a real matrix, non-real eigenvalues of $M$ come in conjugate pairs.
%Hence there is at most one conjugated pair of eigenvalues (and the remaining eigenvalue must be real and positive). %Moreover, as $M$ is a Markov matrix we haveand $\lambda \neq 0$ if $M$ has a real logarithm (Proposition \ref{prop:Culver}).
%
These considerations give rise to Case II and Case III. Any other possibility corresponds to a Markov matrix, with all the eigenvalues real and positive.
Finally, we claim that if the diagonal form is $\dd(1,\lambda,\mu, \mu)$ with $\lambda\neq \mu$, then $\mu\neq 1$. Indeed, if $\mu=1$, then $M-Id$ would have rank 1 because $M$ is diagonalizable. Note that the rows of $M-Id$ vanish, which contradicts the fact that $M-Id$ has no negative entries outside the diagonal.
We are lead to either $\mu=\lambda$ (Case IV) or the three eigenvalues are pairwise distinct (Case I).
%
 %This implies that any diagonalizable $4\times 4$ Markov matrix with a real logarithm lies in one of the cases in Lemma \ref{lemma:DiagCases}.
\end{proof}

Next, we proceed to study the embeddability of  Markov matrices lying in each of the cases in Lemma \ref{lemma:DiagCases}. % in Table \ref{tab:diagonalForms} and Table \ref{tab:diagonalForms2} .

\subsection*{Case I}

\begin{lema}\label{lem_cas1}
Let $M$ be as in Case I with an eigendecomposition $P\dd(1,\lambda_1,\lambda_2,\lambda_3)P^{-1}$  with $\lambda_1,\lambda_2,\lambda_3\in(0,1]$  pairwise distinct and $P\in GL_4(\RR)$. Then $M$ is embeddable if and only if $\Log(M)$ is a rate matrix. Moreover, in this case $\Log(M)$ is the only  Markov generator.
\end{lema}

\begin{proof}
If $\lambda_1,\lambda_2,\lambda_3 \neq 1$, the embeddability of this case is already solved by Corollary \ref{cor:DistinctRealEmbed}. Otherwise, we can assume $\lambda_1=1$ without loss of generality. Under this assumtion, let $Q$ be a Markov generator for $M$.
%It follows from Theorem \ref{thm:CharOfLog} that $Q=P\;A\;\dd(\log_{k_1}(1),\log_{k_2}(\lambda_1),\log_{k_3}(\lambda_2),\log_{k_4}(\lambda_3))\;A^{-1} \;P^{-1}$ for some $k_1,k_2,k_3,k_4\in \ZZ$ and some $A\in Comm^*(\dd(1,\lambda_1,\lambda_2,\lambda_3))$.
 By Remark \ref{rk:IfDifNoCommutator}(i), the eigenvalues of $Q$ are $\log_{k_1}(1),\log_{k_2}(1),\log_{k_3}(\lambda_2),\log_{k_4}(\lambda_3)$ for some $k_i \in \ZZ$.  Since the sum of the rows of $Q$ vanish, $0$ is an eigenvalue of $Q$ and therefore either $k_1=0$ or $k_2=0$. Using that $Q$ is real we deduce that both of them are zero because non-real eigenvalues of $Q$ must appear in conjugate pairs. Again, since $Q$ is real, the eigenvalues of $Q$ corresponding to the non-repeated real eigenvalues of $M$ are their respective principal logarithms, so that $k_3=k_4=0$. As $\Log(M)$ is the only logarithm whose eigenvalues are the principal logarithms of the eigenvalues of $M$ we get $Q=\Log(M)$.
\end{proof}

\subsection*{Case II}\label{sec:NoRepeated4x4}

Markov matrices  $M$ in Case II have non-real eigenvalues and an eigendecomposition as

\begin{equation} \label{eq:Case2}
M=P\; \dd(1,\lambda,\mu,\overline{\mu})\; P^{-1} \text{ with } \lambda \in (0,1], \mu \in \CC \setminus \RR, \ \text{ and } P \in GL_4(\CC).
\end{equation}
{Without loss of generality, we assume $\Imag(\mu)>0$ in order to simplify the notation used in this section.}
If $\lambda \neq 1$,  Proposition \ref{prop:nxnNecessary} claims that the Markov generators of these matrices are of the form
\begin{align*}
 \Ll{k}{M}& = P\; \dd\big(0,\log(\lambda),\log_k(\mu),\overline{\log_k(\mu)} \big) P^{-1}\\
			& =  P \; \dd\big( 0, \log(\lambda), \log(\mu) +2\pi k\ i, \overline{\log(\mu)} -2\pi k\ i  \big) \; P^{-1}.
\end{align*}

% \item $\log^A_k(M)= P\; \dd\big(0,\log(\lambda),\log_k(\mu),\overline{\log_k(\mu)} \big) P^{-1}$ if $Im(\mu)=0$.
%\end{enumerate}

The next result shows that the Markov generators are of this form even if $\lambda=1.$

\begin{prop}\label{Prop_lambda1}
Let $M$ be a Markov matrix with an eigendecomposition $P \; \dd(1,1,\mu,\bar{\mu}) \;P^{-1}$ with $\mu,\bar{\mu} \in \CC$ such that $\mu \neq 0$ and $\Imag(\mu)> 0$. Then,
\begin{itemize}
 \item[(i)] if $\widetilde{P} \; \dd(1,1,\mu,\bar{\mu}) \;\widetilde{P}^{-1}$ is another eigendecomposition of $M$,
 $$P \; \dd(0,0,\log_k(\mu), \overline{ \log_k(\mu)}) \;P^{-1} = \widetilde{P} \; \dd(0,0,\log_k(\mu), \overline{ \log_k(\mu)}) \;\widetilde{P}^{-1} \, ;$$
 \item[(ii)] a matrix $Q$ is a real logarithm of $M$ with rows summing to $0$ if and only if $Q$ has the form $$\Ll{k}{M}=P \; \dd(0,0,\log_k(\mu), \overline{ \log_k(\mu)}) \;P^{-1}.$$
\end{itemize}
\end{prop}

\begin{proof}
$(i)$ If $\widetilde{P} \; \dd(1,1,\mu,\bar{\mu}) \;\widetilde{P}^{-1}$ is another eigendecomposition of $M$, then $\widetilde{P}=PA$ for some matrix $A \in Comm^*(\dd(1,1,\mu,\bar{\mu})).$ As
$$Comm^*(\dd(1,1,\mu,\bar{\mu})) = Comm^*(\dd(0,0,\log_k(\mu), \overline{ \log_k(\mu)})),$$ we obtain the desired result.

$(ii)$ By $(i)$, the definition of $\Ll{k}{M}$ does not depend on $P$ and it is a logarithm of $M$ (see Theorem \ref{thm:CharOfLog}). Note that $(1,1,1,1)^t$ is an eigenvector of $M$ with eigenvalue $1$ because $M$ is a Markov matrix. Hence we can assume that the first column-vector of $P$ is $(1,1,1,1)^t$ and the rows of $\Ll{k}{M}$ sum to $0$.

Conversely, we prove now that any real logarithm $Q$ of $M$ with rows summing to $0$ is of the form $\Ll{k}{M}$.
%Let $Q$ be a real logarithm of $M$ with rows summing to $0$.
From Theorem \ref{thm:CharOfLog} we have that
\begin{eqnarray*}
Q=P\;A\;\dd(\log_{k_1}(1),\log_{k_2}(1),\log_{k_3}(\mu),\log_{k_4}(\bar{\mu}))\;A^{-1} \;P^{-1}
\end{eqnarray*}
 for some $k_1,k_2,k_3,k_4\in \ZZ$ and some $A\in Comm^*(\dd(1,1,\mu,\bar{\mu}))$.  Since the rows of $Q$ sum to $0$ we get $k_1=k_2=0$ as in the proof of Lemma \ref{lem_cas1}.
%Using that $Q$ is real we deduce that $k_2=0$ because otherwise $\log_{k_2}(1)$ would be a \textit{non-real} eigenvalue of $Q$ whose conjugated pair %is not an eigenvalue of $Q$ .
As $Q$ is real, we get that $\log_{k_3}(\mu)$ and $\log_{k_4}(\bar{\mu})$ must be conjugate pairs: $\log_{k_4}(\bar{\mu})=\overline{\log_{k_3}(\mu)}=\log_{-k_3}(\bar{\mu})$ and hence, $k_4=-k_3$. Since the matrix $A$ commutes with $(0,0, \log_{k_3}(\mu), \overline{\log_{k_3}(\mu)})$ (see Remark \ref{rk_comm}), $Q$ is equal to $\Ll{k}{M}$ (taking $k=k_3$).
\end{proof}

Now that we know that all logarithms in Case II are of type $\Ll{k}{M}$, in order to proceed with the study of embedabbility we decompose $\Ll{k}{M}$ as
\begin{equation} \label{eq:L+V}
\Ll{k}{M} = \Log(M) + k \cdot V \text { where }
V= 		P \;\dd(0,0,2\pi i ,-2\pi i)  \; P^{-1}.
\end{equation}

%This description of $\Ll{k}{M}$  is used throughout this section.
Next show that the values of $k$ for which $\Ll{k}{M}$  is a Markov generator form a sequence of consecutive numbers.

\begin{lema}\label{lema:KConnected}
Let $M$ be a Markov matrix as in (\ref{eq:Case2}).  If $\Ll{k_1}{M}$ and $\Ll{k_2}{M}$ are rate matrices with $k_1<k_2$, then $\Ll{k}{M}$ is a rate matrix for all $k \in [k_1,k_2]$.
\end{lema}

\begin{proof}
The proof is immediate because the entries of $\Ll{k}{M} = \Log(M) + k \cdot V $ depend linearly on $k$.
\end{proof}

Note that we could use Lemma \ref{lema:EigenBound} to bound the values of $k$ for which $\Ll{k}{M}$  might be a Markov generator, as we did in Section \ref{sec:NoRepeated}. However, Lemma \ref{lema:KConnected} allows a precise description of those logarithms of $M$  that are Markov generators (not only giving a necessary condition).

\begin{thm} \label{thm:alg}
Let $M$, $P$ and $V$  be as above. Define
%\begin{center}

$$\mathcal{L}:= \displaystyle \max_{(i,j):\ i\neq j,\  V_{i,j} > 0} \left\lceil -\frac{\Log(M)_{i,j}}{V_{i,j}}  \right\rceil,  \qquad \mathcal{U}:= \displaystyle \min_{(i,j):\ i\neq j, \  V_{i,j} < 0} \left\lfloor -\frac{\Log(M)_{i,j}}{V_{i,j}}
  \right\rfloor
$$

and set  $ \ \mathcal{N}:= \{(i,j): i\neq j, \  V_{i,j}=0 \text{ and  } \Log(M)_{i,j}<0\}.$

%\end{center}
 Then,
%\begin{center}
$\Ll{k}{M}$ is a rate matrix if and only if $\mathcal{N} = \emptyset$ and $  \mathcal{L}  \leq k \leq  \mathcal{U}$.
%\end{center}
\end{thm}

\begin{proof}
By (\ref{eq:L+V}) we have  that $\Ll{k}{M}=\mathrm{\Log}(M)+k \cdot V$. Now, assume we choose $k\in \ZZ$ such that $\Ll{k}{M}$ is a rate matrix. In this case, $ \Log(M)_{i,j} + k V_{i,j} \geq 0$ for all $i\neq j$. Hence, for $i\neq j$ we have:

\begin{itemize}

	\item[a)] $0 \leq   \Log(M)_{i,j}$ for all $i,j$ such that $V_{i,j}=0$. In particular $\mathcal{N}=\emptyset$.

	\item[b)] $-\frac{\Log(M)_{i,j}}{V_{i,j}} \leq k$ for all $i,j$ such that $V_{i,j}> 0$. In particular $ \mathcal{L}\leq k$.
	
	\item[c)] $-\frac{\Log(M)_{i,j}}{V_{i,j}} \geq k$  for all $i,j$ such that $V_{i,j}< 0$. In particular $\mathcal{U} \geq k$.

\end{itemize}

Conversely, let us assume that $\mathcal{N} = \emptyset$ and that there is $k\in \ZZ$ such that $ \mathcal{L} \leq k \leq \mathcal{U} $. We want to check that $\Ll{k}{M}$ is a rate matrix. Indeed, take $(i,j)$ with $i\neq j$, then:

%\small
\begin{itemize}
\item[a)] If $V_{i,j} = 0$ we have $\Ll{k}{M}_{i,j}= \Log(M)_{i,j}$. Since $\mathcal{N}=\emptyset$ it follows that $\Log(M)_{i,j}\geq 0$, thus $\Ll{k}{M}_{i,j} \geq 0$.

\item[b)] If $V_{i,j} > 0$, then $\Ll{k}{M}_{i,j}  \geq \Log(M)_{i,j} + \mathcal{L} \cdot V_{i,j} \geq \Log(M)_{i,j} + \frac{-\Log(M)_{i,j}}{V_{i,j}} V_{i,j} = 0.$

\item[c)] If $V_{i,j} < 0$, then $ -\Ll{k}{M}_{i,j} \leq -\Log(M)_{i,j} - \mathcal{U}\cdot V_{i,j} \leq {- \Log(M)_{i,j} - \frac{-\Log(M)_{i,j}}{V_{i,j}} V_{i,j} = 0.}$
Moreover, the rows of $\Ll{k}{M}$ sum to 0, as proved in Prop. \ref{prop:nxnNecessary} and \ref{Prop_lambda1}.
\end{itemize}
%\normalsize
\end{proof}

The theorem above lists all Markov generators of $M$. As an immediate consequence, we get the following characterization of $4\times 4 $ embeddable matrices with a conjugate pair of (non-real) eigenvalues.

\begin{cor}\label{cor:ConjugateEmbeddingCriterion}
Let $M = P\; \dd(1,\lambda,\mu,\bar{\mu})\;P^{-1}$ for some $\lambda\in (0,1]$ and $\mu,\bar{\mu} \in \CC\setminus \RR$. Let $\mathcal{L}$, $\mathcal{U}$ and $\mathcal{N}$ be as in Theorem \ref{thm:alg}. Then, $M$ is embeddable if and only if $\mathcal{N} = \emptyset$ and $\mathcal{L} \leq \mathcal{U}$.
\end{cor}

Now we can prove Theorem \ref{thm_intro} in the introduction using  Lemma \ref{lem_cas1} and Corollary \ref{cor:ConjugateEmbeddingCriterion}:

\vspace*{1mm}

\begin{proofThmintro}
Assume that $M=P\dd(1,\lambda_1,\lambda_2,\lambda_3)P^{-1} $ is a $4\times 4$ Markov matrix with $\lambda_1\in \RR_{>0}$, $\lambda_2\in \CC$, $\lambda_3 \in \CC$ pairwise distinct.
We know that $|\lambda_i|\leq 1$ and, if $M$ is embeddable, $\lambda_i \notin \RR_{\leq 0}$ for any $i=1,2,3$.
Therefore, $M$ lies in Case I if all its eigenvalues are real and in Case II otherwise.

If $M$ lies in Case I, then $M$ is embeddable if and only if $\Log(M)$ is a rate matrix (Lemma \ref{lem_cas1}). As the rows of the principal logarithm of a Markov matrix sum to 0, by setting $V=0$ we have that $\Log(M)$ is a rate matrix  if and only if  $\mathcal{N}=\emptyset$. Moreover, in this case   $\Log(M)$ is the only Markov generator (Lemma \ref{lem_cas1}).

If $M$ lies in Case II, then the statement is precisely Corollary \ref{cor:ConjugateEmbeddingCriterion}. In addition, from Theorem \ref{thm:alg} we obtain that the Markov generators in this case are $\Log_k(M)$ for $k\in [\mathcal{L},\mathcal{U}]$, which coincide with $\Log(M)+2\pi k V$ as defined in the statement of Theorem 1.1.
\end{proofThmintro}

Next, we present an algorithm that solves both the embedding problem and the rate identifiability problem for $4\times 4$ Markov matrices in Cases I and II.
%, which is a dense subset of all $4\times 4$ Markov matrices.
%lying in cases i) and ii) of Proposition \ref{Prop:diagonalForms}. Note that these are a dense subset of all $4\times 4$ Markov matrices.

\begin{rk} \rm
%As claimed in the Introduction,
We already know that the embeddability of a Markov matrix is not always determined by the principal logarithm \cite{K2}. In the $4\times 4$ case, we can prove that the set of embeddable Markov matrices whose principal logarithm is not a Markov generator is not a subset of zero measure; on the contrary, it is a set of full dimension. Moreover, for any $k\in \ZZ$ there is a non-empty Euclidean open set of \emph{embeddable} Markov matrices, all of them in Case II,  whose unique Markov generator is $\Ll{k}{M}$. See \cite{Arxiv}  for details.
%Since the techniques needed to prove this are different to those used so far, the proof is delayed to the Appendix.

%\begin{thm}[Appendix, Theorem \ref{thm:openset}]\label{Thm:Appendix}
%There is a non-empty Euclidean open set of $4\times 4 $ \emph{embeddable} Markov matrices whose principal logarithm is not a rate matrix.
%\end{thm}
\end{rk}

%[Markov generators of $4\times 4$ matrices with a conjugated pair of eigenvalues]

\begin{alg}\label{Alg}\

\begin{algorithm}[H]%\label{Alg}\
\SetKwData{Left}{left}\SetKwData{This}{this}\SetKwData{Up}{up}
\SetKwFunction{Union}{Union}\SetKwFunction{FindCompress}{FindCompress}
\SetKwInOut{Input}{input}\SetKwInOut{Output}{output}

\Input{$M$, a $4\times4$ Markov matrix with distinct eigenvalues as in Thm \ref{thm_intro}.}
\Output{All its Markov generators if $M$ is embeddable, an empty list otherwise.}
\BlankLine

generators=[\ ]\

compute eigenvalues of $M$\

	\If{$M$ has no negative or zero eigenvalue}{

		set $Principal = \Log(M)$\

%\textbf{if} $M$ has a non-real eigenvalue:\

		\If{all the eigenvalues are real}{
			\textbf{add} Principal to generators
		}
		\Else{
			compute $P$ and $V$ as in Thm 1.1\

			compute $\mathcal{L}$, $\mathcal{U}$ and $\mathcal{N}$\

			\If{$\mathcal{N}=\emptyset$}{	
				\For{$k\in \ZZ$ such that $\mathcal{L}\leq k \leq \mathcal{U}$}{
					compute $\Ll{k}{M}=Principal +k \; V$\
		
					\textbf{add} $\Ll{k}{M}$ to generators
				}
			}
		}
	}

\If{generators = [ ]}{\textbf{return} ``M not embeddable''}
\Else{\textbf{return} generators}

%\caption{Markov generators for $4\times4$ matrices with a conjugate pair of eigenvalues.}
\end{algorithm}

\end{alg}

%Note that in Theorem \ref{Thm:Appendix} we saw that there is an open set of embeddable $4\times4$ Markov matrices  whose principal logarithm is not a Markov generator. We complement that result now and postpone the proof for the Appendix:

%\begin{prop}[See Proposition \ref{prop:Allk}] \label{Prop:AllkM}
%For any $k\in \ZZ$ there is a non-empty Euclidean open set of \emph{embeddable} Markov matrices, all of them in case II,  whose unique Markov generator is $\Ll{k}{M}$.
%\end{prop}

\subsection*{Case III} \label{sec:Repeated2}

%Markov matrices in Case III have a real eigenvalue (different than 1) with multiplicity 2, thus they can be written as:
Let $M$ be a Markov matrix as in Case III with an eigendecomposition as
\begin{equation} \label{eq:case3}
 M=P\;\dd(1,\lambda,\mu, \mu)\;P^{-1} \text{ with } \lambda \in (0,1], \  \mu \in[-1,1), \ \mu\neq \lambda,0, \textrm{ and } P \in GL_4(\RR). \end{equation}
{Note  that  the matrix $P$ can be assumed to be real  since all the eigenvalues of $M$ are real}. Note also that this case can be seen as a limit case of Markov matrices with a conjugate pair of complex eigenvalues (case II) and, analogously to that case,  $M$ has infinitely many real logarithms with rows summing to 0. However, in the present case one has to be careful when using Theorem \ref{thm:CharOfLog} in order to take into account the commutant of the diagonal form of $M$.

% has infinite cardinal. %The proofs of the results in this section are included in Appendix \ref{sec:Case3}
%\private{L'altre tamb\'e tenia cardinal infinit, per\`o el pod\'iem treure.}

 %According to Theorem \ref{thm:CharOfLog}, the possible logarithms of $M= P\;\dd(1,\lambda,\mu, \mu)\;P^{-1}$ are $Q_k(A)=(P\; A) \dd(0,\log(\lambda), \log_k(\mu), \overline{\log_k(\mu)}) \; (P \; A)^{-1}$. If we use $\widetilde{P}=P \; A$, the real logarithms of $M$ that diagonalize through $\widetilde{P}$ are of the form $Q_k(A)$.

 % and we could apply the same arguments as in the Case II of Table \ref{tab:diagonalForms} by using $\widetilde{P}$ instead of $P$. However, this is not computationally feasible because the  commutant of the diagonal form of $M$ (where $A$ belongs) has infinite cardinal.

%Note that this corresponds to iterate over the possible eigenvectors of the logarithms of $M$ and find the eigenvalues that make them be rate matrices. What we will do now is the opposite. We will iterate over the possible eigenvalues and find the proper eigenvectors (if they exist). This is possible because the imaginary part of the eigenvectors of any Markov generator of $M$ is bounded (see Lemma \ref{lema:EigenBound}). In particular, there is a finite number of possible eigenvalues of the Markov generators of $M$ and hence the iteration is feasible.

We introduce the following matrices.
 \begin{defi}\label{def:logKRepe}
Let $M,$ $P,$ $\lambda$ and $\mu$ as in (\ref{eq:case3}) above. For any $k\in \ZZ$ and $x,y,z\in \RR$, we define the matrix
\begin{equation*}
 	Q_k(x,y,z)= L + (2\pi k + \Arg{\mu}) \; V(x,y,z), \vspace*{-2mm}
\end{equation*}
where {$  L= P\; \dd(0, \log(\lambda),\log|\mu|,\log|\mu|)\;P^{-1}$} and
$$
 V(x,y,z):= {P \;
 \dd\begin{footnotesize}\left( 0, 0,\begin{pmatrix}
 -y & x\\
 -z & y\\
 \end{pmatrix} \right)
 \end{footnotesize}\;P^{-1}}
.$$
 %and  $ \Arg{\mu}=	\begin{cases}
 % 		0 & \text{if } \mu>0\\
 % 		\pi& \text{if } \mu<0
 % 	\end{cases}$.

%For any $k\in \ZZ$ we define the set $\mathcal{P}_k =\left\lbrace (x,y,z)\in \RR^3: Q_k(x,y,z) \text{ is a rate matrix}\right\rbrace$ and write $\mathcal{V}:=\left\lbrace (x,y,z)\in \RR^3: x>0, \ z>0\text{ and } xz-y^2=1\right\rbrace$.
\end{defi}

\begin{rk}\label{rk:Q0isLog}
\rm If $\mu>0$ we have $Q_0(x,y,z) =\Log(M)$ for all $(x,y,z)\in\RR^3$.
For later use, note that {$V(a x,a y,a z)= a\; V(x,y,z)$ for all $a\in \mathbb{R}$, and hence}
\[Q_k(x,y,z)= \left \{
\begin{array}{ll}
Q_{-k}(-x,-y,-z) & \mbox{ if } \mu>0;  \\
Q_{-k-1}(-x,-y,-z) & \mbox{ if } \mu<0.
\end{array}
\right . \]
\end{rk}

 As in the previous case, we start by enumerating all the real logarithms of $M$ with rows summing to $0$. To this end, we define $\mathcal{V} \subset \RR^3$ as  the algebraic variety
 $$\mathcal{V} =\{(x,y,z) \in \RR^3 \mid \, xz-y^2=1\}.$$
 The next theorem shows that those logarithms with real entries and rows summing to $0$ are of the form $Q_k(x,y,z)$ with $(x,y,z) \in \mathcal{V}$. Furthermore, $\mathcal{V}$ is a $2-$sheet hyperboloid with one of its sheets $\mathcal{V}_-$ in the orthant $x,z<0$ and the other sheet $\mathcal{V}_+$ in the orthant $x,z>0$. The restriction of $(x,y,z)$ to either of these components gives a bijection between  the set of matrices $Q_k(x,y,z)$ and the real logarithms of $Q$ with rows summing to $0$ (other than $\Log(M))$.

\begin{thm}\label{thm:RepeatedLogEnum}\
Let $M$ be a Markov matrix as in  (\ref{eq:case3}). Then, the following are equivalent:\begin{enumerate}[i)]
\item $Q$ is a real logarithm of $M$ with rows summing to $0$;
\item  $Q=Q_k(x,y,z)$ for some $(x,y,z)\in \mathcal{V}$, $k \in \ZZ$.
\end{enumerate}
Moreover, if $Q\neq \Log(M)$ there is a unique $k\in \ZZ$ and a unique $(x,y,z)\in \mathcal{V}_+$ such that $Q=Q_k(x,y,z)$.
\end{thm}

\begin{proof}
i) $\Rightarrow$ ii)
%Since all the eigenvalues of $M$ are real we can always take a real matrix $P$ such that $M= P\; \dd(1,\lambda,\mu,\mu) P^{-1}$.
We know by Theorem \ref{thm:CharOfLog} that any logarithm $Q$ of $M$ is of type $$Q=P\;A\;\dd(\log_{k_1}(1),\log_{k_2}(\lambda),\log_{k_3}(\mu),\log_{k_4}(\mu))\;A^{-1} \;P^{-1}$$ for some $k_1,k_2,k_3,k_4\in \ZZ$ and some $A\in Comm^*(\dd(1,\lambda,\mu,\mu)).$

   Since the rows of $Q$ sum to $0$, $(1,1,1,1)^t$ is an eigenvector of $M$ with eigenvalue $0$. Since non-real eigenvalues of $Q$ must appear in conjugate pairs it follows that $k_1=k_2=0$ (even if $\lambda=1$). Moreover,  we also deduce that  $\log_{k_3}(\mu)$ and $\log_{k_4}(\mu))$ are a conjugate pair. This implies that $k_4= -k_3$ if $\mu >0$ and $k_4= -k_3-1$ if $\mu<0$. Therefore, if we take $k=k_3$, we have
      \small
\begin{equation}\label{eq:LogRepeated}
	\begin{split}
   		   Q & = P\;A\;\dd(\log(1),\log(\lambda),\log_{k}(\mu),\overline{\log_{k}(\mu)})\;A^{-1} \;P^{-1}\\
   		     & = P\;A\;\dd(0,\log(\lambda),\log|\mu|+ (2\pi k +\Arg{\mu} ) i,\log|\mu|- (2\pi k +\Arg{\mu} ) i)\;A^{-1} \;P^{-1}.
   	\end{split}
\end{equation}
\normalsize

If all the eigenvalues of $Q$ are real we deduce that $\Arg{\mu}=0$ and $k=0$. In this case, the eigenvalues of $Q$ are given by the principal logarithm of the respective eigenvalues of $M$ and hence $Q=\Log(M)$.

Now assume that $Q$ has a conjugate pair of complex eigenvalues $\log|\mu| \pm (2\pi k +\Arg{\mu} ) i$. Hence, the third and fourth column-vectors of $P \; A$ must be a conjugate pair (up to scalar product). Furthermore, we have that $P$ is a real matrix and hence it is the third and fourth column-vectors of $ A$ that are a conjugate pair. This fact together with the fact that $A$  commutes with $\dd(1,\lambda,\mu,\mu)$ leads to:

\begin{equation*}
%\begin{small}
\small
A=\begin{pmatrix}
z_1 & 0 & 0 & 0\\
0 & z_2 & 0 & 0\\
0 & 0 & a+bi & z(a-bi)\\
0 & 0 & c+di & z(c-di)\\
\end{pmatrix}
%\end{small}
\end{equation*}
\normalsize
with $z,z_1,z_2 \in \CC\setminus \{0\}$ and $a,b,c,d\in \RR$ satisfying $ad-bc \neq 0$ because $A$ is a non-singular matrix.
We can decompose $A$ as $A=A_1 A_2$ where:
\begin{small}
\begin{equation}\label{eq:Ai}
A_1 =
	\begin{pmatrix}
		1 & 0 & 0 & 0\\
		0 & 1 & 0 & 0\\
		0 & 0 & a & b\\
		0 & 0 & c & d\\
	\end{pmatrix}, \qquad
A_2=\begin{pmatrix}
		z_1 & 0 & 0 & 0\\
		0 & z_2 & 0 & 0\\
		0 & 0 & 1 & z\\
		0 & 0 & i & -zi\\
	\end{pmatrix}.
\end{equation}
\end{small}
Let us define
\begin{align}
J &:=  A_2 \;\dd\left(0,\log(\lambda),\log|\mu|+ (2\pi k +\Arg{\mu} ) i,\log|\mu|- (2\pi k +\Arg{\mu} ) i\right) \; A_2^{-1} \label{eq:Jprod}\\
&=\begin{pmatrix}
		0 & 0 & 0 & 0\\
		0 & \log(\lambda) & 0 & 0\\
		0 & 0 & \log|\mu| & 2\pi k +\Arg{\mu} \\
		0 & 0 & -(2\pi k +\Arg{\mu} ) & \log|\mu|\\
	\end{pmatrix}.\label{eq:Jmat}
\end{align}

Using this notation, the matrix $Q$ in  (\ref{eq:LogRepeated}) can be written as $ Q  =P A_1 J A_1^{-1}P^{-1}$.
Note that $A_1$ commutes with $\dd(0,\log(\lambda),\log|\mu|,\log|\mu|)$ and hence
\small
\begin{equation*}
	Q   = P\;\dd(0,\log(\lambda),\log|\mu|,\log|\mu|) \;P^{-1} + (2\pi k +\Arg{\mu} ) P\;A_1 \; \dd\left(0,0, \begin{pmatrix} 0 & 1 \\ -1 & 0 \end{pmatrix} \right) \;A_1^{-1} \;P^{-1}.
\end{equation*}
\normalsize
A final computation shows that
$A_1 \; \dd\left(0,0, \begin{pmatrix}
0 & 1 \\
-1 & 0
\end{pmatrix} \right) \;A_1^{-1}$ equals $V(x,y,z)$  with
\[x= \frac{a^2+b^2}{ad-bc}, \qquad y=\frac{ac+bd}{ad-bc}, \qquad z=\frac{c^2+d^2}{ad-bc}.\]
It is immediate to show that $xz-y^2=1$, thus $(x,y,z)\in \mathcal{V}$. {This proves that i) implies ii).}
% Moreover, we can take $x,z>0$ because $Q_k(x,y,z)=Q_{-k}(-x,-y,-z)$ if $\mu>0$ and $Q_k(x,y,z)=Q_{-k-1}(-x,-y,-z)$ if $\mu<0$. Hence, we can chose $(x,y,z)$ with $(x,y,z)\in \mathcal{V}_+$. %(analogously $V_-$).

\vspace{2mm}

ii) $\Rightarrow$ i) We know that $Q_k(x,y,z)$ is real by {Definition \ref{def:logKRepe}: it is straightforward to check that $(1,1,1,1)^t$ is an eigenvector with eigenvalue 0 of both $L$ and $V$, and so it also an eigenvector of $Q_k(x,y,z)$ with eigenvalue 0.
% {Moreover, since $M$ is a Markov matrix $(1,1,1,1)^t$ lie in the span of the first column vector of $P$ (or in the span of the first and second columns if $\lambda=1$) and hence the rows of L and $V(x,y,z)$ in Definition \ref{def:logKRepe} sum to zero.
%Therefore, the rows of $Q_k(x,y,z)$ also sum to zero.
}

 Hence it is enough to check  that if $(x,y,z)\in \mathcal{V}$ then $Q_k(x,y,z)$ is a logarithm of $M$. To this end, consider the matrix $J$ introduced in (\ref{eq:Jmat}) and the matrix
$$B := \begin{pmatrix}
1 & 0 & 0 & 0\\
0 & 1 & 0 & 0\\
0 & 0 & 1 & 0\\
0 & 0 & \frac{y}{x} & \frac{1}{x}\\
\end{pmatrix}.$$

If $(x,y,z)\in \mathcal{V}$ then we have $z=\frac{1+y^2}{x}$. A straightforward computation shows that $ P ^{-1} Q_k(x,y,\frac{1+y^2}{x}) P - B J B^{-1} = 0$. Hence, it follows from (\ref{eq:Jprod}) that
\small
\begin{eqnarray*}
Q_k(x,y,\frac{1+y^2}{x})= P\;A\;\dd(0,\log(\lambda),\log|\mu|+ (2\pi k +\Arg{\mu} ) i,\log|\mu|- (2\pi k +\Arg{\mu} ) i)\;A^{-1} \;P^{-1}
\end{eqnarray*}
\normalsize
\noindent with $A= B A_2$ ($A_2$ is defined in (\ref{eq:Ai})).
Since both $B$ and $A_2$ commute with $\dd(1,\lambda,\mu,\mu)$ it follows from Theorem \ref{thm:CharOfLog} that $Q_k(x,y,\frac{1+y^2}{x})$ is a logarithm of $M$, which concludes the first part of the proof.

\vspace{2mm}
In the first part of the proof, we already proved that there exists $k\in \ZZ$ and  $(x,y,z)\in \mathcal{V}$ such that $Q=Q_k(x,y,z)$. By Remark \ref{rk:Q0isLog}, we can take $(x,y,z)\in \mathcal{V}_+$ without loss of generality. To prove that $k$ and $(x,y,z)$ are unique we assume that $Q_k(x,y,z) = Q_{\widetilde{k}}(\widetilde{x},\widetilde{y},\widetilde{z}) $  for some $\widetilde{k}\in \ZZ$ and $(\widetilde{x},\widetilde{y},\widetilde{z})\in \mathcal{V}_+$. In this case, we have $$(2\pi k + \Arg\mu) V(x,y,z)= (2\pi \widetilde{k} + \Arg\mu)V(\widetilde{x},\widetilde{y},\widetilde{z}).$$ Since $Q\neq \Log(M)$ then $(2\pi k + \Arg\mu ) \neq 0$ and hence:
$$x=\frac{2\pi \widetilde{k} + \Arg\mu}{2\pi k + \Arg\mu}\widetilde{x} \qquad y=\frac{2\pi \widetilde{k} + \Arg\mu}{2\pi k + \Arg\mu}\widetilde{y} \qquad z=\frac{2\pi \widetilde{k} + \Arg\mu}{2\pi k + \Arg\mu}\widetilde{z}.$$
  Now, using that $(x,y,z), (\widetilde{x},\widetilde{y},\widetilde{z})\in  \mathcal{V}$ we get $xz-y^2 =  \left( \frac{2\pi \widetilde{k} + \Arg\mu}{2\pi k + \Arg\mu} \right)^2 (\widetilde{x}\widetilde{z}-\widetilde{y}^2)=1$. Moreover, since $x,z, \widetilde{x},\widetilde{z}>0 $, we deduce that $\frac{2\pi \widetilde{k} + \Arg\mu}{2\pi k + \Arg\mu} = 1$, so $\widetilde{k}=k$ and  $(\widetilde{x},\widetilde{y},\widetilde{z})=(x,y,z)$.
\end{proof}

\begin{rk}
\rm Because of Remark \ref{rk:Q0isLog}, every real logarithm of $M$ with rows summing to 0 can also be realized as some $Q_k(x,y,z)$ for a unique $k\in \mathbb{Z}$ and a unique $(x,y,z)\in \mathcal{V}_{-}$. %}{Since $Q_k(x,y,z)=Q_{-k}(-x,-y,-z)$ if $\mu>0$ and $Q_k(x,y,z)=$ $Q_{-k-1}(-x,-y,-z)$ if $\mu<0$, there is also a unique $k\in \ZZ$ and a unique $(x,y,z)\in \mathcal{V}_-$ such that $Q=Q_k(x,y,z)$.}
\end{rk}

In order to characterize those logarithms that are rate matrices, for any $k\in \ZZ$ we define the set $$\mathcal{P}_k =\left\lbrace (x,y,z)\in \RR^3: Q_k(x,y,z) \text{ is a rate matrix}\right\rbrace.$$
Note that the entries of $Q_k(x,y,z)$ depend linearly on $x,y,z$, and hence $\mathcal{P}_k$ is the space of solutions to a system of linear inequalities (i.e. a convex polyhedron). From Theorem \ref{thm:RepeatedLogEnum} we obtain  that the set of Markov generators for a Markov matrix in Case III is $\bigcup_k \mathcal{P}_k \cap \mathcal{V}_+$.
% Although the theorem above does not provide a bound on $k$, we can bound it using Lemma \ref{lema:EigenBound}.
 The following corollary is an immediate consequence of Lemma \ref{lema:EigenBound} and Theorem \ref{thm:RepeatedLogEnum} and shows that there is a finite set of integers $k$ such that $\mathcal{P}_k \cap \mathcal{V}_+ \neq \emptyset$.
%Moreover, excluding the case in Remark \ref{rk:Q0isLog} (that is, $\mu>0$ and $k=0$) we know that $\mathcal{P}_k$ is bounded (See Lemma \ref{lema:PkBounded}).
In Appendix \ref{sec:Case3}  we show a procedure to check whether the intersection $\mathcal{P}_k \cap \mathcal{V}_+$ is not empty and get a point in it.
%Thus, if $\mathcal{P}_k \cap \mathcal{V}_+ \neq \emptyset$ one can easily find a point belonging to it and therefore,

Using the notation introduced in Section \ref{sec:NoRepeated}, if $Q$ is a Markov generator of a Markov matrix $M$ with eigenvalues $1,\lambda$ and $\mu$ as in (\ref{eq:case3}), then it has at most one conjugate pair of non-real eigenvalues, $\log_k(\mu)$ and $\overline{\log_k(\mu)}$. It follows from Lemma \ref{lema:EigenBound} that their imaginary part is bounded by $  \beta_4(\mu) =   -\log|\mu|$  and as consequence, we obtain the next result.
% \tp{He canviat l'antiga $\beta(\mu)$ per la nova $b_4(\log (\mu))$. Caldria repassar que totes les desigualtats i les expressions siguin correctes. Com estem en el cas $n=4$, al denominador del segon membre queda $1=tan(\pi/4)$.}

%
%
%  for any $\mu\in \mathbb{C}$, we write
% %$$b_4(\log(\mu))=\sqrt{2\log(\det(M))\log|\mu|-\log^2|\mu| }= \sqrt{2 \log(\lambda)\log|\mu|+3\log^2|\mu| }.$$
% \begin{eqnarray*}
%  b_4(\log \mu)=\min \{ \sqrt{2\log(\det(M))\log|\mu|-\log^2|\mu| }= \sqrt{2 \log(\lambda)\log|\mu|+3\log^2|\mu|\} }.
% \end{eqnarray*}

\begin{cor} \label{cor:case3RateEnumeration}
Let $M$ be a Markov matrix as in  (\ref{eq:case3}). If  $Q$ is a Markov generator of $M$, then  $Q=Q_k(x,y,z)$ for some $(x,y,z)\in \mathcal{V}_+$ and some $k\in \ZZ$ satisfying  \[\frac{-\Arg \mu + \log|\mu|}{2\pi} \leq k \leq \frac{-\Arg \mu - \log|\mu|}{2\pi}.\]
% $k \in [-\frac{\beta(\mu_j)}{2\pi}, \frac{\beta(\mu_j)}{2\pi}]$ if $\mu>0$ and $k \in [-0.5-\frac{\beta(\mu_j)}{2\pi} , -0.5+ \frac{\beta(\mu_j)}{2\pi}]$ if $\mu<0$.
\end{cor}

As a byproduct, we give an embeddability criterion for $4\times4$ Markov matrices with two repeated eigenvalues.

% \begin{cor}\label{cor:RepeatedEmbedabilitIFF}
% % Let $M$ be a Markov matrix as in  (\ref{eq:case3}) and set $\mathcal{L} = \left\lceil \frac{-Arg(\mu) -b_4(\log(\mu))}{2\pi}  \right\rceil $, $\mathcal{U} = \left\lfloor \frac{-Arg(\mu_j) +b_4(\log(\mu))}{2\pi}  \right\rfloor$. Then, $M$ is embeddable if and only if $\displaystyle\bigcup_{k\in [\mathcal{L}, \mathcal{U}] \cap \ZZ}   \mathcal{P}_k\cap \mathcal{V}_+ \neq \emptyset$.
% \blue{Let $M$ be a Markov matrix as in  (\ref{eq:case3}). If $\mu>0$, define $\mathcal{U}= \lfloor \frac{b_4(\log(\mu))}{2\pi} \rfloor$ and $\mathcal{L}=-\mathcal{U}$; if $\mu < 0$ set  $\mathcal{U} = \lfloor -\frac{1}{2} + \frac{b_4(\log(\mu))}{2\pi} \rfloor$ and $\mathcal{L}=-\mathcal{U}+1$. Then, $M$ is embeddable if and only if $\mathcal{P}_k\cap \mathcal{V}_+ \neq \emptyset$ for some $k\in [\mathcal{L}, \mathcal{U}] \cap \ZZ$.}
% \end{cor}

%\tp{MATEIX COROL.LARI; DIFERENT PRESENTACIO. PER VALORAR}
\begin{cor}\label{cor:RepeatedEmbedabilitIFF}
% Let $M$ be a Markov matrix as in  (\ref{eq:case3}) and set $\mathcal{L} = \left\lceil \frac{-Arg(\mu) -b_4(\log(\mu))}{2\pi}  \right\rceil $, $\mathcal{U} = \left\lfloor \frac{-Arg(\mu_j) +b_4(\log(\mu))}{2\pi}  \right\rfloor$. Then, $M$ is embeddable if and only if $\displaystyle\bigcup_{k\in [\mathcal{L}, \mathcal{U}] \cap \ZZ}   \mathcal{P}_k\cap \mathcal{V}_+ \neq \emptyset$.
Let $M$ be a Markov matrix as in  (\ref{eq:case3}).
\begin{enumerate}
 \item[a)] If $\mu>0$, $M$ is embeddable if and only if $\mathcal{P}_k\cap \mathcal{V}_+ \neq \emptyset$ for some $k$ with
 \[\left \lceil \frac{\log|\mu|}{2\pi} \right \rceil \leq k \leq \left \lfloor \frac{-\log|\mu|}{2\pi} \right \rfloor.\]

  \item[b)] If $\mu<0$, $M$ is embeddable if and only if $\mathcal{P}_k\cap \mathcal{V}_+ \neq \emptyset$ for some $k$ satisfying
 \[ \left \lceil -\frac{1}{2} + \frac{\log|\mu|}{2\pi} \right \rceil \leq k \leq  \left \lfloor -\frac{1}{2} - \frac{\log|\mu|}{2\pi} \right \rfloor.\]

In particular, if $\mu< - e^{-\pi}$ then $M$ is not embeddable.
\end{enumerate}
\end{cor}

\begin{proof}
Since $k\in \ZZ$, the bounds on $k$ are a straightforward consequence of Corollary \ref{cor:case3RateEnumeration}. Indeed, it is enough to take $\Arg \mu = 0$ for $\mu>0$ and $\Arg\mu = \pi$ for $\mu<0$.
In the case of $\mu <0$, it is immediate to check that $ \lceil -\frac{1}{2} + \frac{\log|\mu|}{2\pi} \rceil \ \leq \ \lfloor -\frac{1}{2} - \frac{\log|\mu|}{2\pi} \rfloor$ if and only if $\log|\mu|<- \pi$. Hence, if $\mu< - e^{-\pi}$ there is no $k$ satisfying the embeddability conditions in the statement.
\end{proof}

 \begin{rk}

\rm Example 4.3 in \cite{K2} shows an embeddable Markov matrix as in $(\ref{eq:case3})$ with $\mu =  - e^{-\pi}$. Thus, the bound on Corollary \ref{cor:RepeatedEmbedabilitIFF} is sharp.
 \end{rk}

From Corollary \ref{cor:RepeatedEmbedabilitIFF} we derive an algorithm that tests the embeddability of Markov matrices lying in Case III.\\

\begin{alg}[Markov generators of $4\times 4$ matrices with two repeated eigenvalues]\label{Alg2}\

\begin{algorithm}[H]%\label{Alg2}\

%\caption{Markov generators of $4\times 4$ matrices with two repeated eigenvalues}

\SetKwData{Left}{left}\SetKwData{This}{this}\SetKwData{Up}{up}
\SetKwFunction{Union}{Union}\SetKwFunction{FindCompress}{FindCompress}
\SetKwInOut{Input}{input}\SetKwInOut{Output}{output}

\Input{$M$ (Markov matrix) and $P$ as in (\ref{eq:case3}).}
\Output{One of its Markov generators $Q_k(x,y,z)$ for each $k \in \ZZ$  (if they exist).}
\BlankLine

generators = [ ]\

compute the eigenvalues of $M$: $1$, $\lambda$, $\mu$, $\mu$\

\If{$\det(M) \; >0$ and $\mu \geq -e^{-\pi}$}
{

Compute $L=P\,\dd(0,\log(\lambda),\log|\mu|,\log|\mu|) \; P^{-1}$

	set $\mathcal{L}= \left\lceil \frac{- \Arg{\mu}  + \log|\mu|} {2\pi}   \right\rceil$ and $\mathcal{U}= \left \lfloor \frac{- \Arg{\mu} - |\log(\mu)|   }{2\pi}  \right\rfloor$

	\For{$\mathcal{L}\leq k \leq \mathcal{U}$ $(k\in \ZZ)$:}
	{
		
%\hspace*{10mm}Express $\mathcal{P}_k$ as a system of linear equalities

	\If{ $\mathcal{P}_k \cap \mathcal{V} \neq \emptyset$ (see Appendix \ref{sec:Case3})}
	{
	\textbf{choose} $(x,y,z)\in \mathcal{P}_k \cap \mathcal{V}_+$ (see Appendix \ref{sec:Case3})\

	\textbf{add} $Q_k(x,y,z)=L +k \; V(x,y,z)$ to generators}
	}
}
\If{generators = [ ]}{\textbf{return} ``M not embeddable''}
\textbf{else}{ \textbf{return} generators}
\end{algorithm}

\end{alg}

\begin{rk} \label{rk:Choice}
\rm If $Q_k(x,y,z)\neq \Log(M)$ then each choice of $(x,y,z)\in \mathcal{P}_k \cap \mathcal{V}_+$ in the algorithm above would give a different Markov generator for $M$ (\ref{thm:RepeatedLogEnum}). Thus, the set of \emph{all} Markov generators of $M$ is obtained by considering, for each possible $k$,  all  $(x,y,z)\in \mathcal{P}_k \cap \mathcal{V}_+$ (this can produce infinitely many Markov generators). In Appendix  \ref{sec:Case3} we show how to compute $\# \mathcal{P}_k \cap \mathcal{V}_+$ for a fixed $k$.
\end{rk}

\subsection*{Case IV} \label{sec:Repeated3}

Here, we deal with $4\times 4$ Markov matrices with an eigenvalue of multiplicity $3$ or $4$. This case corresponds to the \textit{equal-input} matrices used in phylogenetics. %The embeddability of this family of matrices is also studied in \cite{Baake}.
{The reader is referred to \cite{Baake} and \cite{BaakeSumner2} for a recent and parallel study on this class of matrices with special emphazis on embeddability.}

\begin{prop} \label{prop:F81Embedd}
Let $M$ be a diagonalizable $4\times4$ Markov matrix with eigenvalues $1,\lambda,\lambda,\lambda$. Then the following are equivalent:
\begin{enumerate}[i)]
	\item $M$ is embeddable.
	\item $\det(M) >0$.
	\item $\Log(M)$ is a rate matrix.
\end{enumerate}
\end{prop}
\begin{proof}

If $M=Id$, that is $\lambda=1$, then  it follows from Theorem \ref{thm:CharOfLog} that $\Log(M)$ is the zero matrix and hence it is a Markov generator for $M$. Moreover, it follows from Corollary \ref{cor:DetBound} the zero matrix is the only Markov generator of the identity matrix.

Now, let us assume $\lambda \neq 1$. Since  $\det(e^Q)=e^{\trace(Q)}$ it follows that  $i) \Rightarrow ii)$. $iii) \Rightarrow i)$ is straightforward, thus to conclude the prove it is enough to check that if $\det(M)>0$ then $\Log(M)$ is a rate matrix.

Since $M$ is a Markov matrix we get that $M-\lambda Id$ is a rank $1$ matrix whose rows sum to $1-\lambda$. Hence:
\begin{equation}\label{eq:F81Matrix}
M = \begin{small}
	\begin{pmatrix}
		a+\lambda  & b & c & d\\	
		a  & b+\lambda & c & d\\
		a  & b & c+\lambda & d\\
		a  & b & c & d+\lambda\\			
	\end{pmatrix} \end{small},
\text{ with } \lambda=1-(a+b+c+d)\in(0,1), \ a,b,c,d \geq 0.
\end{equation}

Let us fix $S\in GL_4(\RR)$ such that $M= S\;\dd(1,\lambda,\lambda,\lambda)\;S^{-1}$. Note that if $\lambda=1$ then $M=Id$ and $\Log(M)=0$ is a rate matrix. On the other hand, if $\lambda\in(0,1)$ then we have:
\begin{equation*}
\begin{split}
 \Log(M) & = S\; \dd(0,\log(\lambda),\log(\lambda),\log(\lambda))\; S^{-1}\\
  		& =  \frac{\log(\lambda)}{\lambda-1} \Big( S\;\dd(1,\lambda,\lambda,\lambda)\;S^{-1} - S\;\dd(1,1,1,1)\;S^{-1} \Big)  =  \frac{\log(\lambda)}{\lambda-1} (M-Id).\\
\end{split}
\end{equation*}
%S\; \begin{tiny}	\left(\frac{x}{\lambda-1},\frac{x\lambda}{\lambda-1} ,\frac{x\lambda}{\lambda-1},\frac{x \lambda}{\lambda-1}\right)\end{tiny}\;S^{-1}  - S\;\left( \frac{x}{\lambda-1},\frac{x}{\lambda-1},\frac{x}{\lambda-1},\frac{x}{\lambda-1}  \right)\;S^{-1}  =
Since $M-Id$ is a rate matrix and $\lambda \in (0,1)$ it follows that $\Log(M)$ is a rate matrix.

\end{proof}

\begin{rk}
\rm In the context of DNA nucleotide-substitution models, these matrices correspond to the Felsenstein81 model \cite{Felsenstein}. The stationary distribution of such matrices is given by $\Pi = (a,b,c,d)/(a+b+c+d)$, where $a,b,c,d$ are as in (\ref{eq:F81Matrix}). When the stationaty distribution is uniform, that is $a=b=c=d$, we recover the Jukes-Cantor model \cite{JC}.
\end{rk}

\subsection{Embeddability of non-diagonalizable $4\times 4$ Markov matrices}

If we restrict the embedding problem to non-diagonalizable $4\times 4$ matrices we have:
\begin{thm} \label{thm:JordanEmbed}
A non-diagonalizable $4\times4$ Markov matrix is embeddable if and only if it has only positive eigenvalues and its principal logarithm is a rate matrix. In this case, it has just one Markov generator.
\end{thm}
\begin{proof}
%First of all, note that in order to the Markov matrix $$M$ has some logarithm it  is necess
The ``if" part is immediate, so we proceed to prove the ``only if" part. Let $M$ be an embeddable non-diagonalizable Markov $4 \times 4$ matrix.  We know that the dominant eigenvalue $1$ has the same algebraic and geometric multiplicity (see \cite[\S 8.4]{Meyer}). Therefore, $M$ has at most one Jordan block of size greater than $1$ and, in this case, its Jordan form is one of the following:
$$\begin{pmatrix}
1 &0 &0  &0 \\
0 & \lambda & 0&0 \\
 0 &0 & \mu & 1 \\
 0 &0 &0 & \mu\\
\end{pmatrix} \text{ with } \mu\neq 0,1,\lambda\neq 0
\qquad \text{ or }\qquad
\begin{pmatrix}
1 &0 & 0 & 0\\
0 & \lambda & 1 &0 \\
0  &0 & \lambda & 1 \\
0  &0 & 0& \lambda\\
\end{pmatrix}\text{ with } \lambda\neq 0,1.$$
{As $M$ is a real matrix, its eigenvalues are necessarily real. Moreover, as $M$ is embeddable,  Proposition \ref{prop:Culver} yields that $\lambda$ and $\mu$ are positive.}
An immediate consequence of Theorem $2$ in \cite{Culver} is that if each Jordan block appears exactly once in its Jordan form, then the only possible real logarithm of $M$ is the principal logarithm.

Hence, if $M$ has a real logarithm other than $\Log(M)$, then the Jordan form of $M$ is $J:=\begin{scriptsize}\begin{pmatrix}
		1 & 0 & 0 & 0\\
		0 & 1 & 0 & 0\\
		0 & 0 & \lambda & 1\\
		0 & 0 & 0 &\lambda\\
	\end{pmatrix} \end{scriptsize}$ with $\lambda \in (0,1)$.

Take $P$ such that $M=P\;J\;P^{-1}$. A more general version of Theorem \ref{thm:CharOfLog} for nondiagonalizable matrices  (see Theorem 1.27  in \cite{Higham}) shows that any logarithm $Q$ of $M$ has the form:
\begin{center}
$Q =  P \; A \; \begin{scriptsize}\begin{pmatrix}
		2\pi k_1 i & 0 & 0 & 0\\
		0 & 2\pi k_2 i & 0 & 0\\
		0 & 0 & \log(\lambda) + 2\pi k_3 i & 1/\lambda\\
		0 & 0 & 0 &\log(\lambda) + 2\pi k_3 i\\
	\end{pmatrix} \end{scriptsize} A^{-1}\; P^{-1}$
\end{center}
 for some $A \in Comm^*(J)$.

It follows that $A$ can be written as $A= \dd(B,Id_2) \dd(c_1,c_2,c_3,c_3)$ with $B\in GL_2(\CC)$. Now, if $Q$ is a rate  matrix it is a real matrix and hence $k_1=-k_2$ and $k_3=0$. Moreover, its rows sum to $0$ and hence $0$ is an eigenvalue of $Q$. Hence, $k_1=-k_2=0$. Thus:
$$Q =  P \; A \; \begin{scriptsize}\begin{pmatrix}
		0  & 0 & 0 & 0\\
		0 & 0 & 0 & 0\\
		0 & 0 & \log(\lambda)  & 1/\lambda\\
		0 & 0 & 0 &\log(\lambda) \\
	\end{pmatrix} \end{scriptsize} A^{-1}\; P^{-1} = P \; \begin{scriptsize}\begin{pmatrix}
		0 & 0 & 0 & 0\\
		0 & 0 & 0 & 0\\
		0 & 0 & \log(\lambda)  & 1/\lambda\\
		0 & 0 & 0 &\log(\lambda) \\
	\end{pmatrix} \end{scriptsize} P^{-1}.$$
and we see that the eigenvalues of $Q$ are the principal logarithms of the eigenvalues of $M$, so that $Q=\Log(M)$.
 \end{proof}

\section{Rate identifiability}\label{sec:RateId}
%!TEX root = main.tex"

Once we know that a Markov matrix arises from a continuous-time model, we want to determine which are its corresponding substitution rates. In other words, given an embeddable matrix we want to know if we can uniquely identify its Markov generator.  Corollary \ref{cor:DetBound} shows that if the determinant of the Markov matrix is big enough, then there is just one generator. However, this is not the case if the determinant is small. Note that a small determinant means that the substitution rates are large or that the substitution process ruled by $M$ has taken a lot of time.

\begin{defi}
An embeddable Markov matrix $M$ has \textit{identifiable rates} if there exists a unique rate matrix $Q$ such that $M=e^Q$.  The \textit{rate identifiability problem} consists on deciding whether a given Markov matrix has identifiable rates or not. 
\end{defi}

\begin{prop}\label{Prop:F81Ident}
Let $M$ be a diagonalizable $4\times4$ embeddable Markov matrix with eigenvalues $1,\lambda,\lambda,\lambda$. If $\det(M)> e^{-6\pi}$, the  rates of $M$ are identifiable and the only generator is $\Log(M)$.
\end{prop}
\begin{proof}
%From Proposition \ref{prop:F81Embedd} we have that $M$ is embeddable if and only if $\lambda>0$. 
Let $Q$ be a Markov generator for $M$. If $\lambda> e^{-2\pi}$ then the real part of the non-zero eigenvalues of $Q$ is greater than $-2\pi$, thus it follows from Lemma \ref{lema:EigenBound} that their imaginary part lies in the interval $(-2\pi,2\pi)$. Since the eigenvalues of $M$ are real and positive this implies that the non-zero eigenvalues of $Q$ are $\log(\lambda)$ and hence $Q= \Log(M)$.
\end{proof}

\begin{rk}\rm
We do not think that this bound is sharp. Up to our knowledge, the largest determinant of a $4\times 4$ embeddable matrix with three repeated eigenvalues and non-identifiable rates is $e^{-12\pi}$, and corresponds to the matrix: :
$$M= \frac{1}{4}
\begin{small}
\begin{pmatrix}
1+3 e^{-4 \pi} & 1-e^{-4 \pi}   & 1-e^{-4 \pi}    & 1-e^{-4 \pi}\\
1-e^{-4 \pi}     & 1+3 e^{-4 \pi} & 1-e^{-4 \pi}    & 1-e^{-4 \pi}\\
1-e^{-4 \pi}     & 1-e^{-4 \pi}   & 1+3 e^{-4 \pi}  & 1-e^{-4 \pi}\\
1-e^{-4 \pi}     & 1-e^{-4 \pi}   & 1-e^{-4 \pi}    & 1+3 e^{-4 \pi}\\
\end{pmatrix}
\end{small}$$

Next we show three Markov generators for it:

$$
\begin{small}
\begin{pmatrix}
-3 \pi & \pi & 2\pi & 0\\
 \pi & -3\pi & 0 & 2\pi\\
 0 & 2\pi & -3\pi & \pi\\
2\pi & 0 & \pi & -3 \pi\\
\end{pmatrix}
\end{small}  \qquad
\begin{small}
\begin{pmatrix}
-3 \pi & \pi & \pi & \pi\\
 \pi & -3\pi & \pi & \pi\\
 \pi & \pi & -3\pi & \pi\\
\pi & \pi & \pi & -3 \pi\\
\end{pmatrix}
\end{small}  \qquad 
\begin{small}
\begin{pmatrix}
-3 \pi & \pi & 0 & 2\pi\\
 \pi & -3\pi & 2\pi & 0\\
 2\pi & 0 & -3\pi & \pi\\
0 & 2\pi & \pi & -3 \pi\\
\end{pmatrix}
\end{small} 
$$
\end{rk}\

 Note that Theorem \ref{thm:nxnCriterion} bounds the number of generators of a Markov matrix with no repeated eigenvalues. Moreover, Algorithm \ref{Alg:nxn} lists all the generators of such a matrix. If we restrict the identifiability problem to $4\times 4$  Markov matrices, we were able to deal with the rate identifiability problem for all the matrices in cases I, II and III, that is, all $4\times4$ matrices except those with an eigenvalue of multiplicity three (Case IV) for which we have Proposition \ref{Prop:F81Ident}. This is summarized in the following table:

% \begin{center}
% \begin{tabular}{|l l l |}
% \hline
% \textbf{Diagonal form of} $\mathbf{M}$ &\textbf{ Embeddability criterion }& \textbf{Number of generators}\\
% 
% \hline
% 
% Case I & $\Log(M)$ is a rate Matrix & One generator.\\
% 
% Case II & $\mathcal{N} = \emptyset$ and $\mathcal{L} \leq \mathcal{U}$ & $U-L+1$ Markov generators\\
% 
% Case III &  $\mathcal{P}_k \cap \mathcal{V} \neq \emptyset$ for some $k\in \ZZ$  &  $\# \displaystyle \cup_{k\in \ZZ} \mathcal{P}_k \cap \mathcal{V}$ generators \\
% 
% Case IV & $\det(M)>0$   & Unknown if $\det(M) \leq e^{-\frac{2 \pi}{\sqrt{3}}}$.\\
% 
% Other diagonal forms & $M$ is not embeddable & \\
% 
% M does not diagonalize &  $\Log(M)$ is a rate Matrix  & One generator.\\
% 
% \hline
% \end{tabular}
% \end{center}

\begin{table}[h]
\begin{center}
\begin{footnotesize}
\begin{tabular}{ccc}
%\hline
%\\

\textbf{Diagonal form of M} & \textbf{Embeddability criterion} & \textbf{Number of generators}\\

\\
\hline
\\

Case I & $\Log(M)$ is a rate Matrix & \jr{1}\\

\\

Case II & $\mathcal{N} = \emptyset$ and $\mathcal{L} \leq \mathcal{U}$ (Cor.  \ref{cor:ConjugateEmbeddingCriterion}) & $\mathcal{U}-\mathcal{L}+1$  (Thm. \ref{thm:alg})\\

\\

Case III &  $\displaystyle \bigcup_k (\mathcal{P}_k \cap \mathcal{V})  \neq \emptyset$  (Cor \ref{cor:RepeatedEmbedabilitIFF}) &  $ \sum_k  \# (\mathcal{P}_k \cap \mathcal{V})$ (Rmk. \ref{rk:Choice})\\

\\

Case IV & $\det(M)>0$  (Prop. \ref{prop:F81Embedd})  & \jr{1} (if $\det(M) > e^{-6 \pi}$)\\

\\

$M$ does not diagonalize &  $\Log(M)$ is a rate Matrix (Thm. \ref{thm:JordanEmbed}) & \jr{1} (Thm. \ref{thm:JordanEmbed})\\

\\

Other diagonal forms & $M$ is not embeddable & $-$\\
\

%\\
%\hline
\end{tabular}
\caption{\label{table_aux} Embeddability test and number of generators for a $4\times 4$ Markov matrix depending on its diagonal form. }
\end{footnotesize}
\end{center}
\end{table}

\section{Discussion}

In this paper we have studied the embeddability and rate identifiability of Markov matrices. 
Our study has lead to a number of results and the development of algorithms that are able to test the embeddability and list the Markov generators of a given Markov matrix, namely Algorithm \ref{Alg:nxn} for any size $n$, and Algorithm \ref{Alg} specifically for $n=4$. 
In this case, the embedding problem has been completely solved. 
It seems natural to think that the next step is to extend the study to $5\times 5$ matrices, at the expense of having to consider two conjugate pairs of complex eigenvalues instead of one and dealing with the difficulties that this causes. 
Most of our results of Section 5 have an immediate generalization to Markov matrices with a single conjugate pair of eiganvalues.

We have used the algorithm \ref{Alg} in a sample of $10^7$ Markov matrices uniformly and independently distributed within the space $\Delta$ of $4\times 4$ Markov matrices. 
%The results are shown in Table \ref{tab:GMvolume}. 
%
As the set of Markov matrices might seem too general for some applied problems (e.g. modeling the substitution of nucleotides in genome), we have previously checked which of the matrices belonged to certain more restrictive families  of matrices that appear in the literature. For instance, for Markov processes on (phylogenetic) trees it is important to restrict to matrices that are \emph{diagonal largest in column} (DLC for short), i.e. Markov matrices whose diagonal entries are the largest entries in each column, see  \cite{Chang}. We denote this set of $4\times 4$ Markov matrices by  $\Delta_{\rm{dlc}}$. Usually, one can even restrict to the set $\Delta_{\rm{dd}}$ of \textit{diagonally-dominant} matrices, that is, matrices $M \in \Delta$ satisfying $M_{ii}\geq 0.5$ for all $i$. If embeddable, these matrices have identifiable rates (\cite{Cuthbert72}). Note that if a matrix is diagonally-dominant, then its off diagonal entries are smaller than or equal to $0.5$. Hence, $\Delta_{\rm{dd}} \subseteq \Delta_{\rm{dlc}}\subseteq \Delta$. 
From a mathematical perspective, it makes sense to target the set of matrices that lie in the connected component of the identity matrix when we remove from $\Delta$ all matrices with determinant equal to $0$. This set is called $\Delta_{\rm{Id}}$, corresponds to matrices with positive determinant, and contains all embedable Markov matrices. 
%It corresponds to the set of Markov matrices with positive eigenvalues if all the transition matrices in $\Delta$ have real eigenvalues. On the contrary, if the model admits transition matrices with non-real eigenvalues, then $\Delta_{\rm{Id}}$  is the set of Markov matrices with positive determinant. This set necessarily includes the multiplicative closure of the transition matrices in the continuous-time version of the model \cite{LMM}.

\begin{table}
  \centering
    \begin{tabular}{|c|c c c |}
    \hline
    
  &  Samples & Embeddable samples  & Percentage of embeddable \\
   \hline
    
   $\Delta$ & $10^7$ & $5774$ & $0.05774$ \\ 
   $\Delta_{\rm{Id}}$ & $4998008$ & $5774$    & $0.11553$   \\ 
   $\Delta_{\rm{dlc}}$  &  $148375$ & $5460$   &  $3.67987$ \\  
   $\Delta_{\rm{dd}}$  &  $2479 $  & $299$ & $12.06132$ \\
   \hline
  \end{tabular}
    \caption{\label{tab:GMvolume} \small Based on a sample of $10^7$ Markov matrices, the first column shows how many sample points lie in each set, the second column shows how many of them are embeddable and the third column displays the corresponding percentage. Embeddability was checked with Algorithm  \ref{Alg}.}
    \end{table}

Therefore, for each matrix generated, we have checked whether it belonged to each of the sets described above ($\Delta_{\rm{dlc}}$, $\Delta_{\rm{dd}}$, $\Delta_{\rm{Id}}$) and, by applying Algorithm \ref{Alg}, we have tested its embeddability. The results are shown in Table \ref{tab:GMvolume}.
To conclude, we would like to notice that as the table shows, the percentage of embeddable Markov matrices is surprisingly small. This result should be taken as a warning signal as it probably has practical consequences related to modeling issues. Indeed, it may warn to reconsider the restriction and use of continuous-time homogeneous models, which seem to be very constrained even for DLC matrices.

\section*{Acknowledgements}
All authors are partially funded by AGAUR Project 2017 SGR-932 and MINECO/FEDER Projects MTM2015-69135, PID2019-103849GB-I00 and MDM-2014-0445. J Roca-Lacostena has received also funding from Secretaria d'Universitats i Recerca de la Generalitat de Catalunya (AGAUR 2018FI\_B\_00947) and European Social Funds.

%The second author started to feel interested in the embedding problem in conversations held with Jeremy Sumner while discussing the definition and properties of the later called Lie Markov models. He wants to express his recognition and gratitude to Jeremy on this matter .

\begin{appendices}

%\section{Embedabble matrices whose principal logarithm is not a rate matrix}\label{sec:Counterexample}
%\input{examples}

\section{Appendix}\label{sec:Case3}

In this appendix, we explain how to find generators for $4\times4$ Markov matrices with two repeated eigenvalues by using Algorithm \ref{Alg2}. More precisely, we explain how to check
whether the intersection $\mathcal{P}_k \cap \mathcal{V}_+$ in Algorithm \ref{Alg2} is empty or not and how to choose a point in it (if not empty).\\

Let $M$ be a diagonalizable $4\times 4$ Markov matrix with a repeated eigenvalue and positive determinant, that is $M=P\;\dd(1,\lambda,\mu,\mu)\; P^{-1}$ for some $P\in GL_4(\RR)$, $\lambda>0$ and $\mu \in [-1,1)$ such that $\mu\neq0$ and $\mu \neq \lambda$. In this case, Theorem \ref{thm:RepeatedLogEnum} yields that each Markov generator other than $\Log(M)$ can be uniquely expressed as $Q_k(x,y,z)$ for some $k\in \ZZ$ and some $(x,y,z)\in \mathcal{P}_k\cap \mathcal{V}_+$.\\

Assume that $k=0$ and $\Arg\mu=0$. Then, $Q_0(x,y,z)$ is equal to the principal logarithm of $M$ for all $(x,y,z)$. Therefore, if the intersection $\mathcal{P}_k \cap \mathcal{V}_+$  is not empty, it is equal to $\mathcal{V}_+$. In this case, the algorithm can choose any point $(x,y,z)\in  \mathcal{V}_+$ such as $(1,0,1)$. For the remainder of this section we assume that $Q_k(x,y,z)\neq \Log(M)$. This assumption is equivalent to assuming that ${2\pi k + \Arg\mu}\neq 0$.\\

% with $\mathcal{V}_+=\{ (x,y,z)\in \RR^3: x,z>0 \text{ and } xz-y^2=1 \}$. %Next we see how to choose $(x,y,z)\in \mathcal{P}_k \cap \mathcal{V}_+$ in Algorithm \ref{Alg2} and study in which cases this intersection is finite or infinite.

 We  denote by $l_{i,j}$ the entries of the matrix $L$ in Definition \ref{def:logKRepe} and by $p_{i,j}$ and $\widetilde{p}_{i,j}$ the entries of $P$ and $P^{-1}$ respectively.  $\mathcal{P}_k$ is the set of solutions to the system of inequalities $Q_k(x,y,z) _{i,j} \geq 0$ for all  $i\neq j$, where $Q_k(x,y,z)= L + (2\pi k +\Arg \mu )V(x,y,z)$. A direct computation shows that the entries of $V(x,y,z)$ depend linearly on $x$, $y$ and $z$:
\begin{equation*}
V(x,y,z) _{i,j}= p_{i,3} \widetilde{p}_{4,j} x -  p_{i,4} \widetilde{p}_{3,j} z +(p_{i,4} \widetilde{p}_{4,j}-p_{i,3} \widetilde{p}_{3,j})y.
\end{equation*}
Hence, the planes $H_{i,j}$ containing the faces of $\mathcal{P}_k$ are given by the equations:
 \begin{equation}\label{eq:Hij}
 p_{i,3} \widetilde{p}_{4,j} x -  p_{i,4} \widetilde{p}_{3,j} z +(p_{i,4} \widetilde{p}_{4,j}-p_{i,3} \widetilde{p}_{3,j})y = \frac{-l_{i,j}}{2\pi k + \Arg{\mu}}.
\end{equation}
From \eqref{eq:Hij} we get that  for each $i\neq j$, the  faces of two polyhedra $\mathcal{P}_{k_1}$ and $\mathcal{P}_{k_2}$ ($k_1,k_2 \in \ZZ$) corresponding to the $(i,j)$-entry of $Q_k(x,y,z)$ are necessarily parallel.\\ %Thus if $P_k$ is unbounded for some positive $k$, then it is unbounded for all positive $k$. The same applies for negative values. Moreover, depending on the sign of $\mu$ we have that $(x,y,z)\in\mathcal{P}_k$ if and only if $(-x,-y,-z)\in \mathcal{P}_{-k}$ or $(-x,-y,-z)\in \mathcal{P}_{-k-1}$. Hence, if $P_k$ is unbounded for some $k$ (other than $k=0$ if $\mu>0$) then it is unbounded for all other $k$. In particular, there are values of $k$ such that $k> -\frac{\sqrt{3} \log(\det(Q)}{4 \pi}$ and $\mathcal{P}_k  \neq \emptyset$, which is in contradiction with Lemma \ref{lema:EigenBound}.

%Next we proceed, to show how to find points in $\mathcal{P}_k\cap \mathcal{V}_+$ for any $k\in \ZZ$.

%Fixed $k$, each point in $\mathcal{P}_k$ satisfying $f(x,y,z)=0$ corresponds to a Markov generator $Q_k(x,y,z)$.  By  modifying the algorithm above one can obtain a parametrization of all the Markov generators of $M$ with parameters $k\in [\mathcal{L},\mathcal{U}]\cap \ZZ $ and $(x,y,z)\in \mathcal{P}_k \cap \mathcal{V}$.To chose a value of $(x,y,z)$ we proceed as follows:

Let us  define $f(x,y,z)= xz-y^2-1$ so that $\mathcal{V}=\{(x,y,z)\in \RR^3 \mid f(x,y,z)=0\}$. Note that $(x,y,z)\in \mathcal{V}_+$ if and only if $f(x,y,z)=0$ and $x,z>0$.   Next we show how to find points in $\mathcal{P}_k\cap \mathcal{V}_+$. To do so, we  evaluate $f(x,y,z)$ at the vertices, on edges and on faces of $\mathcal{P}_k$ according to the following procedure:
 \begin{itemize}
  \item[Step 1] Evaluate $f(x,y,z)$ at each of the vertices of $\mathcal{P}_k$.

$-$ If there is a pair of vertices $v_1$ and $v_2$ such that $f(v_1) f(v_2)<0$, then $\mathcal{V}_{+}$ cuts $\mathcal{P}_k$ and hence there are infinitely many generators. To find one of them, we restrict $f$ to the line defined by $v_1$ and $v_2$, and find a point $P=(x,y,z)$ in the segment between $v_1$ and $v_2$ such that $f(P)=0$.

$-$ If there is not such a pair of vertices but there is some vertex $v$ satisfying $f(v)=0$ then we take $P=v$. %In this case, if there is exactly one vertex satisfying $f(v)=0$ and $\mathcal{V}_{+} \cap (\mathcal{P}_k\setminus {v}) \neq \emptyset$ in a neighbourhood of $v$ then there are infinitely many of generators (for the current value of $k$). Otherwise,  the number of generators (for the current value of $k$) is the number vertices such that $f(v)=0$.

$-$ If the evaluation of $f$ at all the vertices of $\mathcal{P}_k$ has the same sign (and none is equal to $0$), then we proceed to Step 2.

\item[Step 2] Find the vanishing points of $f(x,y,z)$ on the edges of $\mathcal{P}_k$.  To do so, we restrict $f$ to the lines containing these edges and look for solutions $P=(x,y,z)$ of $f=0$ lying in the corresponding edge.

$-$ If we find two solutions in such a line, $\mathcal{V}_{+}$ cuts the interior of $\mathcal{P}_k$ and hence there are infinitely many generators.

% If all the solutions lying in the edges have multiplicity 2, then we can choose the values of $(x,y,z)$ corresponding to any of these points.

%In this case, if $\mathcal{V}_{+}$ intersects the edges of $\mathcal{P}_k$ in exactly one point $P=(x,y,z)$ and the intersection of $\mathcal{V}_{+}$ and $\mathcal{P}$ restricted to a neighbourhood of $P$ is exactly $P$, then $Q_k(P)$ is the only Markov generator for the current value of $k$ (otherwise, there are infinitely many generators).

$-$ If the edges of $\mathcal{P}_k$ do not intersect $\mathcal{V}_{+}$, then we proceed to Step 3.

\item[Step 3] We assume that $\mathcal{V}$ does not intersect any edge of $\mathcal{P}_k$.  For $i\neq j$ consider the intersection $H_{i,j}\cap \mathcal{V}_+$, where $H_{i,j}$ is the plane defined by \eqref{eq:Hij}.

If this intersection is not empty, choose a point in it (see next paragraph) and check whether it belongs to $\mathcal{P}_k$ or not.
%Note that in the previous step we got that the edges of $\mathcal{P}_k$ do not intersect $\{f(x,y,z)=0\}$.
This intersection lies either completely in the corresponding face of the polyhedron or completely outside the polyhedron.

  $-$ If we find a point $P=(x,y,z)$ which belongs to the polyhedron in this way, then $M$ has infinitely many generators and $Q_k(x,y,z)$ is one of them.

  $-$ If we fail to find a point in any of the faces, then $M$ has no generator with the current value of $k$.
 \end{itemize}
To conclude, we give some insight on how to find $P \in H_{i,j}\cap \mathcal{V}_+$ (when this intersection is not empty). For ease of reading we write  $Ax +B y + C z = D$ for the equation of $H_{i,j}$ (see \eqref{eq:Hij} for the precise coefficients). Given $(x,y,z)\in \mathcal{V}$, we can write $z=\frac{1+y^2}{x}$ because $x\neq 0$. Therefore,  by multiplying the equation $Ax+By+Cz=D$ by $x$ and rearranging the terms in the equality, we conclude that
\small
$$ C  y^2 + (B x)  y + (Ax^2- D x+ C) =0 {\normalsize \text{ if and only if }} \left(x,y,\frac{1+y^2}{x}\right)\in \mathcal{V}\cap H_{i,j}  .$$
\normalsize
Hence, $\mathcal{V_+}\cap H_{i,j}$ is not empty if there exists $x>0$ for which the discriminant $$\Disc := (B^2-4 A C) \; x^2 + (4 \; C D)\;x- 4 C^2$$ is non-negative.  We study below whether this is possible depending on the coefficients of $\Disc$:
 \begin{enumerate}[i)]

\item \label{it:Inter1} If $B^2-4 A C<0$, compute the roots of $\Disc=0$. If they are non-real or negative then  $\mathcal{V}_+ \cap H_{i,j} = \emptyset$. If both roots are real and  positive, then all values between them satisfy $\Disc \geq 0$.   Note that in this case, both roots have necessarily the same sign (because $\mathcal{V}\cap \{x=0\} = \emptyset$).

\item \label{it:Inter2} If $B^2-4 A C>0$, then $\Disc > 0$ when
$x\gg 0$.
%$x\rightarrow +\infty$.

%Note that,  since $\mathcal{P}_k$ is bounded then this implies that the intersection of $\mathcal{V}_+$ with the face $(i,j)$ of the polyhedron is empty.

\item \label{it:Inter3} If $B^2-4 A C=0$ and $C  D>0 $, then $\Disc > 0$ when $x\gg 0$. %$x\rightarrow +\infty$.
% Since $\mathcal{P}_k$ is bounded then this implies that the intersection of $\mathcal{V}_+$ with the face $(i,j)$ of the polyhedron is empty.

\item  If $B^2-4 A C=0$, $C  D\leq 0 $ and $C\neq0$, then $\Disc<0$ for all $x\geq0$. In this case, we have $\mathcal{V}_+ \cap H_{i,j} = \emptyset$.

\item  \label{it:Inter5} If $\Disc$ is identically $0$, then $H_{i,j}\cap \mathcal{V}_+$ is unbounded with respect to $x$, that is, for any $x_0>0$, there are points $P=(x_0,y,z)\in H_{i,j}\cap \mathcal{V}_+$.

%Hence, since $\mathcal{P}_k$ is bounded, this implies that the intersection of $\mathcal{V}_+$ with the face $(i,j)$ of the polyhedron is empty.

\end{enumerate}

Note that in cases \ref{it:Inter2}), \ref{it:Inter3}) and \ref{it:Inter5}) we have that $\Disc \geq 0$ when $x\gg 0$ and hence the $H_{i,j} \cap \mathcal{V}_+$ is unbounded with respect to $x$. Therefore, if we know that the polyhedron $\mathcal{P}_k$ is bounded, then the intersection of $\mathcal{V}_+$ with the face $(i,j)$ is necessarily empty in any of these cases.

\end{appendices}

\bibliographystyle{alpha}
\bibliography{biblio}

\end{document}